\documentclass[11pt,letterpaper]{amsart}

\usepackage[utf8]{inputenc}
\usepackage{enumerate}%
\usepackage{enumitem}
\usepackage{amsmath,amssymb}
\usepackage{cite}%%%
\usepackage{comment}
\usepackage{amsfonts}%%
\usepackage{mathrsfs}%%
\usepackage{xcolor}
\usepackage[colorlinks, linkcolor = blue, anchorcolor = blue, citecolor = blue]{hyperref}
\usepackage{hyperref}
\usepackage{cleveref}
\usepackage{needspace}

\theoremstyle{definition}
\newtheorem{definition}{Definition}[section]
\newtheorem{remark}[definition]{Remark}

\newtheorem{theorem}[definition]{Theorem}

\newtheorem{lemma}[definition]{Lemma}
\newtheorem{proposition}[definition]{Proposition}
\newtheorem{corollary}[definition]{Corollary}

\DeclareMathOperator{\cof}{cof}
\DeclareMathOperator{\tr}{tr}
\DeclareMathOperator{\supp}{supp}
\DeclareMathOperator{\diag}{diag}
\DeclareMathOperator{\diver}{div}

\theoremstyle{plain}

\newcommand{\R}{\mathbb R}

\newcommand{\Z}{\mathbb Z}
\newcommand{\T}{\mathbb T}

\newcommand{\loc}{\mathrm{loc}}
\newcommand{\A}{\mathcal A}
\newcommand{\XX}{\mathcal X}
\newcommand{\YY}{\mathcal Y}

\newcommand{\e}{\mathrm e}
\newcommand{\ii}{\mathrm i}
\newcommand{\norm}[1]{\left\lVert#1\right\rVert}

\numberwithin{equation}{section}

\newcommand{\beq}{\begin{equation}}
\newcommand{\eeq}{\end{equation}}
\newcommand{\beqs}{\begin{eqnarray*}}
\newcommand{\eeqs}{\end{eqnarray*}}
\newcommand{\beqn}{\begin{eqnarray}}
\newcommand{\eeqn}{\end{eqnarray}}
\newcommand{\beqa}{\begin{array}}
\newcommand{\eeqa}{\end{array}}

\title[Unique continuation]{Unique continuation and its failure for the Monge--Amp\`ere equation}

\author{Xiaorui Lu}
\address{School of Mathematical Sciences, Zhejiang University, Hangzhou 310058, China}
\email{12435035@zju.edu.cn}

\date{}

\begin{document}

\begin{abstract}
We study unique continuation for Alexandrov solutions with common
Monge--Amp\`ere density bounded above and below. In the plane,
we prove that agreement on an open set forces global equality. For every bounded domain \(\Omega\subset\mathbb R^n\) with \(n\ge3\),
we construct two distinct uniformly convex \(C^{1,1}(\Omega)\) solutions that agree on a nonempty open subset of \(\Omega\).
\end{abstract}
\maketitle

%%%%%%%%%%%%%%%%%%%%%%%%%%%%%%%%%%%%%%%%%%%%%%%%%%%%%%%%
\baselineskip16pt
\parskip3pt

%%%%%%%%%%%%%%%%%%%%%%%%%%%%%%%%%%%%%%%%%%%%%%%%%%%%%%%
%%%%%%%%%%%%%%%%%%%%%%%%%%%%%%%%%%%%%%%%%%%%%%%%%%%%%%%

% \tableofcontents
\section{Introduction}

Mooney posed the following unique continuation problem for bounded densities in his survey \cite[\S6.1]{Mooney2018}.
Let \(u,v\) be convex functions on a connected open set
\(\Omega\subset\mathbb R^n\) with
\begin{equation}\label{eq:question}
  M_u=M_v=f\,dX \qquad\text{in }\Omega.
\end{equation}
If \(0<\lambda\le f\le\Lambda<\infty\) almost everywhere and \(u=v\)
on a nonempty open subset of \(\Omega\), must they coincide throughout
\(\Omega\)? Here \(M_w(E)=|\partial w(E)|\) is the Monge--Amp\`ere measure
of \(w\), the Lebesgue measure of its subgradient image.

For smooth Monge--Amp\`ere equations, Aronszajn proved
uniqueness for \(C^{4,1}\) elliptic solutions with matching Cauchy data
on a portion of \(C^{4,1}\) boundary \cite{Aronszajn1957}.
Mooney proved unique continuation in every dimension for Alexandrov
solutions with strictly positive \(C^{1,\alpha}\) density \cite{Mooney2015}. For measurable
\(f\) with \(\lambda\le f\le\Lambda\), Jian and Tu proved that an ordered pair \(u\ge v\)
coincides if $u, v$ touch at an interior point of strict convexity
of either one \cite{JianTu2025}. Mooney's question allows \(u-v\) to
change sign.

De Philippis, Guerra and Tione used quantitative unique continuation 
to recover planar interior $W^{2,1+\delta}$ estimates \cite{DePhilippisGuerraTione2026}.
Uhlmann and Zimmermann established uniqueness for nonnegative $L^1$ densities 
from generalized boundary measurements on bounded uniformly convex $C^{1,1}$ domains \cite{uhlmann2026recoveryroughsourcesinverse}.

Equality of the Monge-Ampère measures \eqref{eq:question} cannot be
replaced by equality of the Hessian determinants almost
everywhere.
Let $u(X)=|X|^2/2$ and $v(X)=u(X)+\max\{X_1,0\}$.
Both \(u\) and \(v\) have Hessian determinant \(1\) almost everywhere and agree on \(\{X_1<0\}\).
\(v\) has additional mass on the interface, since its subgradient
image of \(\{0\}\times E\) is \([0,1]\times E\). Related examples arise
in the comparison of Brenier and Alexandrov solutions
\cite{CortesPoza2026,DePhilippisFigalli2014}.

In one dimension, \(M_w=w''\) in distributions. Equation \eqref{eq:question} implies that
\(u-v\) is affine, so agreement on an interval suffices even when the
common measure is singular. In the plane, the density bounds can be local.

\begin{theorem}
\label{thm:lowdim}
Let $\Omega\subset\mathbb{R}^2$ be a connected open set,
and $f$ be a measurable function on $\Omega$.
Assume that for every open set $U\Subset\Omega$ there exist
constants $0<\lambda_U\le\Lambda_U<\infty$ such that
\[
    \lambda_U\le f(x)\le\Lambda_U
    \quad\text{for a.e. }x\in U.
\]
Suppose that $u,v:\Omega\to\mathbb{R}$ are Alexandrov
solutions of \eqref{eq:question}.
If $u=v$ on a nonempty open subset of $\Omega$,
then $u=v$ throughout $\Omega$.
\end{theorem}

Planar strict convexity\cite{Alexandrov1942} and interior \(C^{1,\alpha}\) regularity
\cite{Caffarelli1991} permit a local partial Legendre
transform in \(x\), which was used in Lewy and Heinz's works on Monge-Amp\`ere equations in dimension two
and can be traced back to Darboux in the 19th century \cite{TrudingerWang2008, Liu2021Interior}. For \(C^{1,1}\) solutions
in a class with \(w_{xx}\ge c>0\), Guan derived a divergence equation for
the inverse coordinate \cite[\S5, Lemma~18]{Guan1997}. The difference \(h=u^*-v^*\) of the transformed potentials satisfies
\[
  \partial_p(a\,\partial_p h)+\partial_y^2 h=0,\qquad
  a\,\partial_p h=\int_{X_v(p,y)}^{X_u(p,y)} f(s,y)\,ds .
\]
Here \(\partial_p h=X_u-X_v\), so \(a\) is the average of \(f(\cdot,y)\)
between the two inverse coordinates. 
% We take \(a\) equal to the local lower bound for \(f\) if \(X_u=X_v\). 
The coefficient matrix \(\operatorname{diag}(a,1)\)
is measurable and uniformly elliptic, so the planar theorem of Alessandrini
and Magnanini propagates equality in the original variables
\cite{Alessandrini1992,AM1994}.

For \(n\ge3\), the partial Legendre transform leaves a nonlinear determinant in
the \(n-1\) transverse variables. In these dimensions unique continuation
fails even for solutions with Hessians arbitrarily close to the identity.
Write \(B_1^n\) for the unit ball,
\(B_{1,+}^n=B_1^n\cap\{X_1>0\}\), and \(Q_n(X)=|X|^2/2\).

\begin{theorem}
\label{thm:main}
For each $n\ge3$ and $0<\varepsilon<1$, there are
$u_n,v_n\in C^{1,1}(B_1^n)$ and a measurable $f_n$ such that
\begin{enumerate}
\item $(1-\varepsilon/4)I_n\le D^2u_n,D^2v_n\le
 (1+\varepsilon/4)I_n$ almost everywhere;
\item $u_n=v_n=Q_n$ on $B_{1,+}^n$, while
 $u_n(t,0,\ldots,0)>v_n(t,0,\ldots,0)$ for $-1<t<0$;
\item for every Borel set $E\subset B_1^n$,
\begin{equation}\label{eq:main-measure}
 M_{u_n}(E)=M_{v_n}(E)=\int_E f_n\,dX;
\end{equation}
\item $1-\varepsilon\le f_n\le1+\varepsilon$ almost everywhere,
 and $f_n=1$ on $B_{1,+}^n$.
\end{enumerate}
\end{theorem}

For \(f\equiv c>0\), Mooney's theorem gives unique continuation in every
dimension \cite[Theorem~1.3]{Mooney2015}. For \(n\ge3\), every
\(L^\infty(B_1^n)\)-neighborhood of \(c\) contains a density for which
continuation fails. The pair \(c^{1/n}u_n,c^{1/n}v_n\) has common
density \(g=c f_n\), with \(\|g-c\|_{L^\infty(B_1^n)}\le c\varepsilon\)
and \(g=c\) a.e. on \(B_{1,+}^n\).

\begin{corollary}\label{cor:interval}
Let $\Omega\subset\mathbb{R}^n$ be a bounded domain with $n\ge 3$.
For any $0<m<M<\infty$, there exists a measurable density $f$ on
$\Omega$, with $m\le f\le M$ almost everywhere, such that \eqref{eq:question}
in $\Omega$ admits two distinct uniformly convex Alexandrov solutions
$u,v\in C^{1,1}(\Omega)$ that agree on a nonempty open subset of
$\Omega$.
\end{corollary}

This gives a negative answer to Mooney's unique continuation question
\cite[\S6.1, item~2]{Mooney2018} in every dimension $n\ge3$ on every bounded domain.
Pli\'s constructed a three-dimensional linear elliptic counterexample
in nondivergence form \cite{Plis1963}. Miller and Mandache obtained
divergence-form counterexamples with H\"older coefficients
\cite{Mandache1998,Miller1973,Miller1974}. We adapt Miller's oscillations
of increasing tangential frequency. For \(u_3=\phi+w\) and \(v_3=\phi-w\),
the determinants agree precisely when
\[
  \operatorname{cof}(D^2\phi):D^2w+\det D^2w=0 .
\]
We correct both \(\phi\) and \(w\) to solve this identity. The area formula for the injective Lipschitz gradients gives
\eqref{eq:main-measure}. The difference \(h=u_n-v_n\) is a nonzero
solution of \(\operatorname{div}(A\nabla h)=0\) vanishing on \(B_{1,+}^n\),
where \(A=\int_0^1\operatorname{cof}D^2(v_n+s h)\,ds\). This uniformly
elliptic matrix is divergence-free and can be made
arbitrarily close to \(I_n\) in \(L^\infty\).

Section~\ref{pre} recalls some useful results
used in Section~\ref{sec:po2} to prove Theorem~\ref{thm:lowdim}. Sections~\ref{sec:algebra}--\ref{sec:nonlinear} build
oscillating functions in three variables and use a right inverse for
the linear equation to solve the nonlinear determinant identity.
Section~\ref{sec:completion} extends the solutions across the interface, verifies their
Monge--Amp\`ere measures and lifts the examples to every dimension \(n>3\).

\section{Preliminaries}
\label{pre}

\subsection{Monge--Amp\`ere measures and regularity}

We first recall the precise definition of Alexandrov solutions. Any convex
function \(w : \Omega \subset \mathbb{R}^n \to \mathbb{R}\) has an associated
Borel measure \(M_w\), called the Monge--Amp\`ere measure, defined by
\(M_w(E) = |\partial w(E)|\),
where \(|\partial w(E)|\) represents the Lebesgue measure of the image of
the subgradients of \(w\) in \(E\), see \cite{Gutierrez2001}. We say \(p \in \mathbb{R}^n\)
is a subgradient of \(w\) at \(x\) if it is the slope of some supporting
hyperplane to the graph of \(w\) at \(x\). If \(w \in C^2\), then
\(|\partial w(E)| = \int_E \det D^2 w\,dx.\)
Given a Borel measure \(\mu\), we say that \(w\) is an Alexandrov solution to
\(\det D^2 w = \mu\) if \(M_w = \mu\).

\begin{lemma}
\label{pl:lem:regularity}
If $M_w=g\,dx\,dy$ with $0<\lambda\le g\le\Lambda$ on a
planar domain $\Omega$, then the Alexandrov solution $w$ is locally strictly convex and
$w\in C^{1,\alpha}_{\loc}$ for some $\alpha>0$.
If $g$ is smooth on $\Omega$, then $w$ is smooth there and
$D^2w$ is positive definite.
\end{lemma}

Local strict convexity follows from Alexandrov planar theorem
\cite{Alexandrov1942}; see also
\cite[Remark~2.19]{DePhilippisFigalli2014} for a modern formulation
and historical attribution.
The local $C^{1,\alpha}$ regularity follows from Caffarelli \cite{Caffarelli1991}.
When $g$ is smooth, Caffarelli's interior $C^{2,\alpha}$ regularity
theorem \cite{Caffarelli1990Regularity}, followed by standard
elliptic bootstrapping, gives smoothness of $w$.
Positive definiteness of $D^2w$ then follows from convexity and
$\det D^2w=g>0$.
% We use these results without reproducing their proofs.

\begin{lemma}
\label{pl:lem:approximation}
Let $\overline B\subset\Omega$ be a ball and suppose that
$M_w=g\,dx\,dy$ near $\overline B$, with the bounds above.
There exist smooth $g_j$ and convex
$w_j\in C(\overline B)\cap C^\infty(B)$ such that
\begin{gather*}
 \lambda\le g_j\le\Lambda,\quad g_j\longrightarrow g\text{ in }L^1(B),
 \quad \det D^2w_j=g_j,\quad w_j=w\text{ on }\partial B,\\
 w_j\longrightarrow w\text{ uniformly on }\overline B,
 \qquad Dw_j\longrightarrow Dw\text{ locally uniformly in }B.
\end{gather*}
\end{lemma}

Here $g_j$ may be obtained by mollification with a nonnegative
kernel after extending $g$ by $\lambda$ outside $B$.
The classical existence theory for the Alexandrov Dirichlet
problem goes back to Bakelman \cite{Bakelman1957} and
Alexandrov \cite{Alexandrov1958}.
For the existence and stability results used here, we refer to
\cite[Theorem~2.14 and Proposition~2.16]{Figalli2017}.
Lemma~\ref{pl:lem:regularity} supplies interior smoothness.
Since $w\in C^1(B)$, local uniform convergence of the gradients
follows from the classical convergence theorem for differentiable
convex functions; see \cite[Theorem~25.7]{Rockafellar1970}.

For context, we recall the established interior Sobolev regularity
$w\in W^{2,1}_{\loc}$ for strictly convex solutions with bounded
positive density \cite{DePhilippisFigalli2013}, together with the
almost everywhere identity $\det D^2w=g$, see \cite[Lemma~2.1]{TrudingerWang2005}.
The stronger $W^{2,1+\eta}$ estimates 
\cite{DePhilippisFigalliSavin2013,Schmidt2013} are not needed below.

\subsection{Partial Legendre transform}

Let \(w=w(x,y)\) be a strictly convex \(C^1\) function. Its partial Legendre transform with respect to \(x\) is
\[
w^*(p,y)=\sup_x\{px-w(x,y)\},
\]
with \(p=w_x(x,y)\). Write \(x=X_w(p,y)\), then
\[
w^*(p,y)=pX_w(p,y)-w(X_w(p,y),y).
\]
The transformed function $w^*$ is \(C^1\) with
\begin{equation}\label{pl:eq:firstderivatives}
      (w^*)_p=X_w,\qquad
(w^*)_y=-w_y(X_w(p,y),y).
\end{equation}
In particular,
\[
x=(w^*)_p(p,y),\qquad
w(x,y)=xp-w^*(p,y).
\]
In addition, if \(w\in C^2\) and \(w_{xx}>0\), then
\begin{equation}\label{eq:plt}
\begin{gathered}
  (X_w)_p = (w^*)_{pp} = \frac{1}{w_{xx}}, \qquad
  (X_w)_y = (w^*)_{py} = -\frac{w_{xy}}{w_{xx}}, \\[4pt]
  (w^*)_{yy} = -w_{yy} + \frac{w_{xy}^2}{w_{xx}}
             = -\frac{\det D^2 w}{w_{xx}}.
\end{gathered}
\end{equation}
For Alexandrov solutions, these second-order calculations are performed on smooth approximants.
See \cite{Heinz1956,Heinz1959} and
\cite[\S5.1.4]{Figalli2017} for further information.

Fix $R=(a,b)\times(c,d)$ with $\overline R\subset B\Subset\Omega$.
For a strictly convex $C^1$ function $w$ near $\overline R$, define
\begin{equation}\label{pl:eq:coordinate}
 T_w(x,y)=(w_x(x,y),y),\qquad
 O_w=\{(p,y):c<y<d,\ w_x(a,y)<p<w_x(b,y)\}.
\end{equation}
Strict monotonicity of $w_x(x,y)$ in $x$-coordinate makes $T_w$ a
bijection from $R$ onto $O_w$. Write its inverse as
$(X_w(p,y),y)$ and set
\begin{equation}\label{eq:w-star}
 w^*(p,y)=\max_{x\in[a,b]}\{px-w(x,y)\}
         =pX_w(p,y)-w(X_w(p,y),y).
\end{equation}

\begin{lemma}
\label{pl:lem:coordinates}
The set $O_w$ is open, $T_w:R\to O_w$ is a homeomorphism.
If strictly convex $C^1$ functions $w_j$ converge to $w$ in
$C^1(\overline R)$, then every $K\Subset O_w$ is contained in
$O_{w_j}$ for all sufficiently large $j$, and, uniformly on $K$,
\begin{equation*}
 X_{w_j}\longrightarrow X_w,\qquad w_j^*\longrightarrow w^*,
 \qquad Dw_j^*\longrightarrow Dw^*.
\end{equation*}
\end{lemma}

    For each fixed $y$, the function $x\mapsto w_x(x,y)$ is continuous
and strictly increasing. Thus $T_w$ is a continuous bijection from
$R$ onto $O_w$, and $O_w$ is open by \eqref{pl:eq:coordinate}.
% The maximum in \eqref{pl:eq:plt} is attained at the unique point $X_w(p,y)$. 
% Berge's maximum theorem
% \cite[Chap.~VI, \S~3]{Berge1959} shows that $X_w$ is continuous.
Since \(w_x(x,y)\) is continuous and strictly increasing in \(x\), the inverse coordinate \(X_w(p,y)\) is continuous in \((p,y)\).
Consequently,
\(T_w^{-1}(p,y)=(X_w(p,y),y)\)
is continuous, so $T_w$ is a homeomorphism.
Since $T_{w_j}\to T_w$ locally uniformly on $R$,
\cite[Lemma~2.5]{Smith2006} gives both $K\subset O_{w_j}$ for all
sufficiently large $j$ and $X_{w_j}\to X_w$ uniformly on $K$.
The uniform convergence of $w_j^*$ follows from 
\cite[Remark~2.3(2)]{DaskalopoulosSavin2009}.
Since $X_{w_j}\to X_w$ uniformly and
$w_j\to w$ in $C^1(\overline R)$, \eqref{pl:eq:firstderivatives} and
the uniform continuity of $w_y$ give $Dw_j^*\to Dw^*$ uniformly on $K$.

\subsection{Planar unique continuation}
\begin{lemma}[Alessandrini--Magnanini]
\label{pl:lem:planar-ucp}
Let $Q\subset\R^2$ be a disk and let $A$ be a bounded measurable,
symmetric, uniformly positive definite matrix field. A weak solution
$h\in W^{1,2}_{\loc}(Q)$ of $\diver(A\nabla h)=0$ that vanishes on a
nonempty open subset is identically zero. The same conclusion holds if
its zero set has positive planar measure.
\end{lemma}

The published paper
\cite[Theorem~2.1]{AM1994} writes $h=H\circ\chi$ locally, with $H$
harmonic and $\chi$ quasiconformal; the corresponding strong
continuation statement is \cite[Corollary~2.2]{AM1994}.
The conclusion for positive-measure zero set follows directly
from \cite[Lemma~2.3]{HarrachPohjolaSalo2019}.
The 1992 report \cite{Alessandrini1992} predates the published proof;
this chronology is explicitly documented by the author in
\cite[Introduction]{Alessandrini2012}.

\subsection{Area formula and determinant identities}
For a uniformly convex $u\in C^{1,1}(\Omega)$ on a convex domain,
the Lipschitz gradient is injective. The area formula therefore
gives, for every Borel set $E$,
\begin{equation}\label{eq:area}
 M_u(E)=|\nabla u(E)|=\int_E\det D^2u\,dX.
\end{equation}
We use the precise area formula in
\cite[Theorem~3.2.3]{Federer1969}. In particular, there is no singular interface mass for
the uniformly convex $C^{1,1}$ functions constructed below.

For symmetric $3\times3$ matrices we shall use the elementary
determinant identities
\begin{equation}\label{eq:odd-det}
    \det(A+B)-\det(A-B)=2(\cof A:B+\det B),
\end{equation}
\begin{equation}\label{eq:cofactor}
    \cof(I_3+H)=I_3+(\tr H)I_3-H+\cof H.
\end{equation}
The Piola identity gives $\diver(\cof D^2q)=0$ distributionally for
$q\in W^{2,\infty}_{\loc}$; it is the null-Lagrangian identity 
applied to $\nabla q$, with the Sobolev formulation
obtained by approximation \cite{BallCurrieOlver1981}.

\section{Proof of Theorem~\ref{thm:lowdim}}
\label{sec:po2}

Fix a rectangle $R=(a,b)\times(c,d)$ and a ball $B$ such that
$\overline R\subset B\Subset\Omega$.
Choose constants $0<\lambda\le\Lambda<\infty$ such that
$\lambda\le f\le\Lambda$ almost everywhere on a neighborhood of
$\overline B$.
We choose a Borel representative of $f$ taking values in
$[\lambda,\Lambda]$ there.
For each solution $w$, we use the same notation as Section~\ref{pre}.

Suppose that $w,f\in C^\infty(B)$, $D^2w>0$, and
$\det D^2w=f$.
By \eqref{eq:plt}, 
\begin{equation*}\label{pl:eq:transformed-potential}
 w^*_{yy}+f(X_w,y)w^*_{pp}=0.
\end{equation*}
Differentiating this identity with respect to $p$, and using
$w^*_p=X_w$, gives
\begin{equation*}\label{pl:eq:guan}
 (X_w)_{yy}
 +\partial_p\bigl(f(X_w,y)(X_w)_p\bigr)=0.
\end{equation*}
This is the weak divergence-form equation appearing in
\cite[Lemma~18, equation~(71)]{Guan1997}.

\begin{proposition}
\label{pl:prop:weak-identity}
Suppose that $u$ and $v$ are Alexandrov solutions of
\[
 M_u=M_v=f\,dx\,dy
\]
on $B$.
Let $Q\Subset O_u\cap O_v$ be a disk and set
\(h=u^*-v^*.\)
Define
\begin{equation}\label{pl:eq:coefficient}
 a(p,y)=
 \begin{cases}
 \displaystyle
 \frac{1}{X_u(p,y)-X_v(p,y)}
 \int_{X_v(p,y)}^{X_u(p,y)}f(s,y)\,ds,
 &X_u(p,y)\ne X_v(p,y),\\[8pt]
 \lambda,
 &X_u(p,y)=X_v(p,y).
 \end{cases}
\end{equation}
Then $a$ is measurable with
$ \lambda\le a\le\Lambda$,
and
\begin{equation}\label{pl:eq:difference}
 \partial_p(a h_p)+h_{yy}=0
 \qquad\text{in }\mathcal D'(Q).
\end{equation}
Equivalently, for every $\eta\in C_c^\infty(Q)$,
\begin{equation}\label{pl:eq:weak-test}
 \int_Q\bigl(a h_p\eta_p+h_y\eta_y\bigr)\,dp\,dy=0.
\end{equation}
\end{proposition}

\begin{proof}
First suppose that $f,u,v\in C^\infty(B)$,
$D^2u,D^2v>0$, and $\det D^2u=\det D^2v=f$.
Since $h_p=X_u-X_v$, the definition of $a$ gives
\[
 a h_p=\int_{X_v}^{X_u}f(s,y)\,ds.
\]
Differentiating this identity with respect to $p$ yields
\[
 \partial_p(a h_p)
 =f(X_u,y)(X_u)_p-f(X_v,y)(X_v)_p
 =-u^*_{yy}+v^*_{yy}
 =-h_{yy}.
\]
This proves \eqref{pl:eq:difference} in the smooth case.
Multiplication by $\eta\in C_c^\infty(Q)$ and integration by parts
give \eqref{pl:eq:weak-test}.

For general Alexandrov solutions $u,v$, the coefficient $a$ is
measurable, since it can be written as
\[
 a(p,y)=\int_0^1
 f\bigl(X_v(p,y)+t(X_u(p,y)-X_v(p,y)),y\bigr)\,dt
\]
on $\{h_p\ne0\}$,
while $a=\lambda$ on $\{h_p=0\}$.
Because $f$ is Borel measurable with values in $[\lambda,\Lambda]$,
we have $\lambda\le a\le\Lambda$.
Moreover, \(a h_p=\int_{X_v}^{X_u}f(s,y)\,ds\)
also holds where $h_p=0$.

Take the Dirichlet approximants $f_j,u_j,v_j$ from
Lemma~\ref{pl:lem:approximation}, using the same density $f_j$
for both solutions, and set
\[
 h_j=u_j^*-v_j^*,\qquad
 X_{u_j}=(u_j^*)_p,\qquad
 X_{v_j}=(v_j^*)_p.
\]
For fixed $\eta\in C_c^\infty(Q)$,
Lemma~\ref{pl:lem:coordinates} gives
$\overline Q\subset O_{u_j}\cap O_{v_j}$ for large $j$ and
\[
 u_j^*\longrightarrow u^*,\qquad
 v_j^*\longrightarrow v^*
 \quad\text{in }C^1(\overline Q).
\]
Define $a_j$ by \eqref{pl:eq:coefficient}, with
$f,X_u,X_v$ replaced by $f_j,X_{u_j},X_{v_j}$.
Then
\[
 a_j(h_j)_p=\int_{X_{v_j}}^{X_{u_j}}f_j(s,y)\,ds.
\]
The smooth identity yields
\begin{equation}\label{pl:eq:smooth-test}
 \int_Q\bigl(a_j(h_j)_p\eta_p+(h_j)_y\eta_y\bigr)
 \,dp\,dy=0.
\end{equation}
The bound $f\le\Lambda$ gives
\[
 \begin{split}
 &|a_j(h_j)_p-a h_p|
 =\left|
 \int_{X_{v_j}}^{X_{u_j}}f_j(s,y)\,ds
 -
 \int_{X_v}^{X_u}f(s,y)\,ds
 \right|\\
 \le &
 \int_a^b|f_j(s,y)-f(s,y)|\,ds +\Lambda|X_{u_j}-X_u|
       +\Lambda|X_{v_j}-X_v|.
 \end{split}
\]
Hence,
\begin{equation}\label{pl:eq:composition-convergence}
 \begin{split}
 &\|a_j(h_j)_p-a h_p\|_{L^1(Q)}\\
 \le &
 \operatorname{diam}(Q)\|f_j-f\|_{L^1(R)}+
 \Lambda|Q|
 \left(
 \|X_{u_j}-X_u\|_{L^\infty(Q)}
 +
 \|X_{v_j}-X_v\|_{L^\infty(Q)}
 \right)
 \longrightarrow0.
 \end{split}
\end{equation}
Lemma~\ref{pl:lem:coordinates} also gives
\(\|(h_j)_y-h_y\|_{L^\infty(Q)}\longrightarrow0,\)
so
\[
 \left|
 \int_Q\bigl((h_j)_y-h_y\bigr)\eta_y\,dp\,dy
 \right|
 \le
 \|(h_j)_y-h_y\|_{L^\infty(Q)}
 \|\eta_y\|_{L^1(Q)}
 \longrightarrow0.
\]
Together with \eqref{pl:eq:composition-convergence}, this justifies
passing to the limit in \eqref{pl:eq:smooth-test} and proves
\eqref{pl:eq:weak-test}, hence \eqref{pl:eq:difference}.
\end{proof}

Only the smooth approximants are differentiated. In particular, this
argument does not require a chain rule for $f(X_w(p,y),y)$ with
merely measurable $f$. Guan weak equation
\cite[Lemma~18]{Guan1997} is an antecedent under different regularity
assumptions.

\begin{proof}[Proof of Theorem~\ref{thm:lowdim}]
By Lemma~\ref{pl:lem:regularity}, $u$ and $v$ are
locally strictly convex and belong to $C^1_{\loc}(\Omega)$. Set
\[
 E=\{z\in\Omega:u=v\text{ in a neighborhood of }z\}.
\]
This set is nonempty and open. We prove that it is relatively closed.
Take $z_0=(x_0,y_0)\in\overline E\cap\Omega$ and $z_j\in E$ with
$z_j\to z_0$. Equality of $u,v$ and their gradients on $E$, followed
by continuity, gives
\begin{equation*}
 u(z_0)=v(z_0),\qquad Du(z_0)=Dv(z_0).
\end{equation*}
Select a rectangle $R$ and a ball $B$ with
$z_0\in R$, $\overline R\subset B$, and $\overline B\subset\Omega$,
and choose positive bounds $\lambda,\Lambda$ for $f$ on a
neighborhood of $\overline B$. Then
\[
 T_u(z_0)=(u_x(z_0),y_0)=(v_x(z_0),y_0)=T_v(z_0)
\]
lies in $O_u\cap O_v$. Choose a disk $Q$ about this point with
$\overline Q\subset O_u\cap O_v$.

For all sufficiently large $j$, we have $z_j\in E\cap R$ and
$T_u(z_j)\in Q$. Thus $T_u(E\cap R)\cap Q$ is nonempty; it is open
since $T_u$ is a homeomorphism.
For $(p,y)\in T_u(E\cap R)\cap Q$, let
$(x,y)=T_u^{-1}(p,y)$. Since $u=v$ and $u_x=v_x$ at $(x,y)$,
we have $T_v(x,y)=(p,y)$ and $T_v^{-1}(p,y)=(x,y)$.
Formula~\eqref{eq:w-star} therefore gives
\[
 u^*(p,y)=px-u(x,y)=px-v(x,y)=v^*(p,y).
\]
Thus $h=u^*-v^*$ vanishes on the nonempty open set
$T_u(E\cap R)\cap Q$.
Since $h\in C^1(Q)\subset W^{1,2}_{\loc}(Q)$ and
$\diag(a,1)$ is bounded, measurable, symmetric, and uniformly elliptic,
Lemma~\ref{pl:lem:planar-ucp} applies to $h$, implying $h=0$, hence $u^*=v^*$, throughout $Q$.
Their $p$ derivatives also agree on $Q$, so
\eqref{pl:eq:firstderivatives} yields
\[
 T_u^{-1}(p,y)=(u^*_p(p,y),y)
             =(v^*_p(p,y),y)=T_v^{-1}(p,y).
\]
The inverse image $T_u^{-1}(Q)=T_v^{-1}(Q)$ is therefore an
open neighborhood of $z_0$. For any $(x,y)$ in this neighborhood,
take $(p,y)=T_u(x,y)=T_v(x,y)$, then
\[
 u(x,y)=xp-u^*(p,y)=xp-v^*(p,y)=v(x,y).
\]
It follows that $z_0\in E$. Therefore, $E$ is nonempty, open, and
relatively closed in the connected set $\Omega$, so $E=\Omega$.
\end{proof}

\section{Construction in dimension three}\label{sec:algebra}

We first construct two functions in three variables, which are discussed in Sections~\ref{sec:algebra}--\ref{sec:nonlinear}. The extension to
all $n\ge3$ is included in the proof of Theorem~\ref{thm:main} in
Section~\ref{sec:completion}. The analytic statement we need is the following.

\begin{proposition}\label{prop:three}
For every $\delta>0$, there are $\phi,w\in C^{1,1}(B_1^3)$ such that
\begin{equation}\label{eq:small-hessians}
    \|D^2\phi-I_3\|_{L^\infty(B_1^3)}
 +\|D^2w\|_{L^\infty(B_1^3)}\le\delta,
\end{equation}
\begin{equation}\label{eq:boxed}
      \cof(D^2\phi):D^2w+\det D^2w=0
 \quad\hbox{a.e. in }B_1^3.
\end{equation}
Moreover, $\phi=Q_3$ and $w=0$ on $B_{1,+}^3$, while
$w(t,0,0)>0$ for $-1<t<0$.
\end{proposition}

Write $X=(t,x,y)$. We seek
\begin{equation*}\label{eq:pair}
 u_3=Q_3+\psi+w,\qquad v_3=Q_3+\psi-w.
\end{equation*}
In Proposition~\ref{prop:three}, this means
$\phi=Q_3+\psi$; the construction below uses $\psi$ directly.
The determinant identity \eqref{eq:odd-det} shows that the required
equality is equivalent to
\begin{equation}\label{eq:exact}
 \cof(I_3+D^2\psi):D^2w+\det D^2w=0.
\end{equation}

The proof has four parts. We choose a function $ah_0$ whose spatial
frequency increases while amplitude decays, with Laplacian
relatively small. We then write $w=a(h_0+z)$ and construct a right
inverse for the linear equation in $(\psi,z)$, allowing us to solve \eqref{eq:exact} by a contraction map. Finally, we
extend $\psi,w$ across $t=0$ and check convexity and Monge--Amp\`ere
measures. The use of increasing frequencies is motivated by Miller\cite{Miller1973, Miller1974}. The linear inverse and nonlinear estimates needed
for \eqref{eq:exact} are proved below.

All estimates for the construction are carried out in three dimensions before the lift in Section~\ref{sec:completion}. Constants
denoted by $C$ may change from line to line. Unless otherwise specified, they depend
only on fixed functions $S, \zeta, \chi$. They
are uniform in $N$, the interval indices, and all Fourier frequencies. Hessian bounds without a Fourier norm use the operator norm.
Throughout the construction, $D$ differentiates in all three variables,
whereas $\Delta_{x,y}=\partial_{xx}+\partial_{yy}$.
All corrections will be $2\pi$-periodic in $x,y$. We work on
$\T^2=(\R/2\pi\Z)^2$ and then regard these functions as periodic
functions on $\R^2$, allowing us to solve an ordinary differential
equation for each spatial frequency.

\subsection{Fourier norms and derivative estimates}\label{sec:norms}
For an interval $I\subset\R$ and a Fourier series
\[
 b(t,x,y)=\sum_{m\in\Z^2}b_m(t)\e^{\ii(m_1x+m_2y)},
\]
define
\begin{equation}\label{eq:A-norm}
 \norm b_{\A(I)}=\sum_{m\in\Z^2}\norm{b_m}_{L^\infty(I)}.
\end{equation}
If $b$ is real-valued, for almost every $t$,
\( \overline{b_m(t)}
 =b_{-m}(t).\)
Moreover, $\norm b_{L^\infty(I\times\T^2)}\le\norm b_{\A(I)}$.
For a vector $v=(v_i)$ and a matrix $M=(M_{ij})$, we define
\[
 \norm v_{\A(I)}=\sum_i\norm{v_i}_{\A(I)},\qquad
 \norm M_{\A(I)}=\sum_{i,j}\norm{M_{ij}}_{\A(I)}.
\]
For $b,c\in\mathcal A(I)$, absolute convergence gives $(bc)_m=\sum_{q\in\mathbb Z^2}b_qc_{m-q}$. Hence, by Tonelli theorem,
\begin{equation}\label{eq:algebra}
\|bc\|_{\mathcal A(I)}\le\sum_{m,q\in\mathbb Z^2}\|b_q\|_{L^\infty(I)}\|c_{m-q}\|_{L^\infty(I)}=\|b\|_{\mathcal A(I)}\|c\|_{\mathcal A(I)}.
\end{equation}
For matrices, with $A:B=\sum_{i,j}A_{ij}B_{ij}$, this gives
\begin{equation*}\label{eq:matrix-algebra}
 \norm{A:B}_{\A(I)}\le\norm A_{\A(I)}\norm B_{\A(I)},\qquad
 \norm{\tr A}_{\A(I)}\le\norm A_{\A(I)}.
\end{equation*}

For $k\ge1$, define the second-order norm
\begin{equation}\label{eq:second-norm}
 \norm b_{2,k;I}=\sum_m\left[
 (k^2+|m|^2)\norm{b_m}_{\infty,I}
 +(k+|m|)\norm{b_m'}_{\infty,I}
 +\norm{b_m''}_{\infty,I}\right].
\end{equation}
Time derivatives here are weak derivatives. This norm controls the
function, its gradient, and its Hessian at the scale $k$,
\begin{equation}\label{eq:derivative-control}
 k^2\norm b_{\A(I)}+k\norm{Db}_{\A(I)}
       +\norm{D^2b}_{\A(I)}\le C\norm b_{2,k;I}.
\end{equation}
This estimate follows directly from the Fourier differentiation rules.
For each coefficient,
\[
 \begin{array}{ll}
 (\partial_t b)_m=b_m',&(\partial_{x_\ell}b)_m=\ii m_\ell b_m,\\
 (\partial_{tt}b)_m=b_m'',&
 (\partial_{t x_\ell}b)_m=\ii m_\ell b_m',\qquad
 (\partial_{x_\ell x_r}b)_m=-m_\ell m_r b_m,
 \end{array}
\]
where $x_1=x,x_2=y$. Set $s_m=|m_1|+|m_2|$. The contribution of frequency $m$ to the left side of
\eqref{eq:derivative-control} is exactly
\( (k^2+ks_m+s_m^2)\norm{b_m}_{\infty,I}
 +(k+2s_m)\norm{b_m'}_{\infty,I}
 +\norm{b_m''}_{\infty,I}.\)
Since $s_m\le\sqrt2|m|$ and $ks_m\le(k^2+s_m^2)/2$, we have
\( k^2+ks_m+s_m^2\le3(k^2+|m|^2), 
 k+2s_m\le3(k+|m|).\)
Summing over $m$ proves \eqref{eq:derivative-control}, with a constant
independent of $k$.

The finite norm in \eqref{eq:second-norm} also gives a $C^{1,1}_{\mathrm{loc}}$ representative.
Indeed, \eqref{eq:second-norm} implies that \(b\in W^{2,\infty}_{\mathrm{loc}}\), thus both \(b\) and its first weak
derivatives admit locally Lipschitz representatives; see
\cite[Theorem~3.38, p.~71]{Hunter2014}.
Consequently, on every compact convex set $K\Subset I\times\mathbb R^2$, this representative satisfies
\[
 |Db(Y)-Db(X)|\le C\|b\|_{2,k;I}|Y-X|,
 \qquad \text{for} \,X,Y\in K.
\]
We use this representative throughout.

We will multiply solutions by smooth functions of time. The formulas
\[
 (\zeta b_m)'=\zeta b_m'+\zeta'b_m,\qquad
 (\zeta b_m)''=\zeta b_m''+2\zeta'b_m'+\zeta''b_m
\]
give
\begin{equation}\label{eq:cutoff-norm}
 \norm{\zeta b}_{2,k;I}
 \le C\left(\norm\zeta_\infty+k^{-1}\norm{\zeta'}_\infty
                    +k^{-2}\norm{\zeta''}_\infty\right)
                    \norm b_{2,k;I}.
\end{equation}

\subsection{Estimates for matrices}
Let $H,B,\widetilde H,\widetilde B$ be
$3\times3$ matrix functions with entries in $\A(I)$.
All norms $||\cdot||$ in this subsection are $\A(I)$ norms.
Since cofactor is a difference of two products and determinant
is a sum of six signed products,
\begin{equation}\label{eq:polynomial-size}
 \norm{\cof H}\le C\norm H^2,\qquad
 \norm{\det B}\le C\norm B^3.
\end{equation}

To estimate differences, apply
\(ab-\widetilde a\widetilde b
=(a-\widetilde a)b+\widetilde a(b-\widetilde b)\)
to the two products in each cofactor.
\eqref{eq:algebra}, followed by summation, gives
\begin{equation}\label{eq:cofactor-difference}
 \|\operatorname{cof}H-\operatorname{cof}\widetilde H\|
 \le C\bigl(\|H\|+\|\widetilde H\|\bigr)
       \|H-\widetilde H\|.
\end{equation}
Similarly, apply
\(abc-\widetilde a\widetilde b\widetilde c
=(a-\widetilde a)bc
+\widetilde a(b-\widetilde b)c
+\widetilde a\widetilde b(c-\widetilde c)\)
to the six products in the determinant. Applying \eqref{eq:algebra} and summing therefore yields
\begin{equation}\label{eq:det-difference}
 \|\det B-\det\widetilde B\|
 \le C\bigl(\|B\|+\|\widetilde B\|\bigr)^2
       \|B-\widetilde B\|.
\end{equation}
These estimates do not require any derivatives of the matrix entries.

The spaces with norms \eqref{eq:A-norm}
and \eqref{eq:second-norm} are complete.
Let \(\|\cdot\|_*\) denote either norm in \eqref{eq:A-norm}
and \eqref{eq:second-norm}, with \(k\ge1\) fixed.
If \((b^{(n)})\) is Cauchy, completeness of the coefficient spaces
gives, for every \(m\),
\(b_m^{(n)}\to b_m\) in \(L^\infty(I)\) or
\(W^{2,\infty}(I)\), respectively.
Let \(b\) have Fourier coefficients \(b_m\). Fatou lemma gives
\[
 \|b\|_*\le\liminf_{n\to\infty}\|b^{(n)}\|_*<\infty,\qquad
 \|b^{(n)}-b\|_*
 \le\liminf_{\ell\to\infty}
       \|b^{(n)}-b^{(\ell)}\|_*
 \longrightarrow0,
\]
when $n\to\infty$. Thus both spaces are complete.
If \(b^{(n)}_{-m}=\overline{b^{(n)}_m}\), then
\(b_{-m}=\lim_n b^{(n)}_{-m}
=\overline{\lim_n b^{(n)}_m}=\overline{b_m}\);
hence their real subspaces are closed.

\subsection{The oscillating functions}\label{sec:waves}

We now construct $h_0$ and $a$. The function $h_0$ alternates between
oscillations in $x$ and in $y$.
The amplitude \(a(t)\) tends to zero as \(t\uparrow0\) ensuring \(a h_0\to0\) uniformly
in \(x,y\).
We choose specific frequencies and transition intervals so that
\(a^{-1}\Delta(a h_0)\) is relatively small; see \eqref{eq:base-estimates}.

Fix a smooth nondecreasing function $S:\R\to[0,1]$ equal to zero
on $(-\infty,0]$ and to one on $[1,\infty)$. Let $N$ be a large positive integer and set
\begin{equation*}\label{eq:frequencies}
 k_j=(N^4+j)^2,\qquad
 t_j=-2N^2\sum_{r=j}^\infty\frac1{k_r},\qquad j\ge1.
\end{equation*}
The function $s\mapsto(N^4+s)^{-2}$ is decreasing, so 
\begin{equation*}\label{eq:time-bounds}
 \frac{2N^2}{N^4+j}\le -t_j
 \le\frac{2N^2}{N^4+j-1}.
\end{equation*}
In particular, $t_j\uparrow0$, $t_{j+1}-t_j=2N^2/k_j$, and
$-1/4<t_1<0$ when $N \ge 3$. Adjacent frequencies satisfy
\begin{equation}\label{eq:ratio}
 \frac{k_{j+1}}{k_j}
 =1+\frac2{N^4+j}+\frac1{(N^4+j)^2}
 \le1+3N^{-4}.
\end{equation}

Let $e_j=(1,0)$ for odd $j$ and $e_j=(0,1)$ for even $j$. On
$[t_j,t_{j+1}]$, define
\begin{equation*}
    \theta_j(t)=\frac\pi2 S\!\left(\frac{k_j(t-t_j)}{N^2}\right),
\end{equation*}
\begin{equation}\label{eq:h0}
    h_0(t,x,y)=\cos\theta_j(t)\cos(k_je_j\cdot(x,y))+\sin\theta_j(t)\cos(k_{j+1}e_{j+1}\cdot(x,y)),
\end{equation}
\begin{equation*}\label{eq:kappa}
    \kappa(t)=k_j+(k_{j+1}-k_j)S\!\left(\frac{k_j(t-t_j)}{N^2}\right).
\end{equation*}
For $t<t_1$, set $h_0=\cos(k_1x)$ and $\kappa=k_1$.
The change from one frequency to the next takes place in
$[t_j,t_j+N^2/k_j]$. During
$[t_j+N^2/k_j,t_{j+1}]$, the formulas reduce to
\[
 h_0=\cos(k_{j+1}e_{j+1}\cdot(x,y)),\qquad \kappa=k_{j+1}.
\]
These are exactly the values at the start of the next interval.
Since \(S\in C^\infty(\mathbb R)\), all positive-order
\(t\)-derivatives of \(h_0\) and \(\kappa\) vanish at each
junction \(t_j\), so the piecewise formulas define smooth functions for \(t<0\).

Define
\begin{equation*}\label{eq:amplitude}
 g(t)=\int_{-1}^t\kappa(s)\,ds,
 \qquad a(t)=N^{-2}k_1^{-2}\e^{-g(t)},\qquad \text{when } \,t<0.
\end{equation*}
We will look for
\begin{equation*}\label{eq:w-form}
w(t,x,y)=a(t)\bigl(h_0(t,x,y)+z(t,x,y)\bigr).
\end{equation*}
Only $z$ and $\psi$ remain to be found. Differentiating $a$ gives
$a'=-\kappa a$, and $a''=(\kappa^2-\kappa')a$.
Define the two normalized Hessians
\begin{equation*}\label{eq:B-def}
 B_0=a^{-1}D^2(ah_0),\qquad B_z=a^{-1}D^2(az).
\end{equation*}
The product rule gives every entry of the second matrix:
\begin{equation}\label{eq:B-entries}
 \begin{split}
 (B_z)_{tt}&=z_{tt}-2\kappa z_t+(\kappa^2-\kappa')z,\\
 (B_z)_{t\ell}&=z_{t\ell}-\kappa z_\ell
       \quad(\ell=x,y),\\
 (B_z)_{\ell r}&=z_{\ell r}\quad(\ell,r=x,y).
 \end{split}
\end{equation}
In particular,
\(\tr B_z=z_{tt}-2\kappa z_t+(\kappa^2-\kappa')z+z_{xx}+z_{yy}.\)

\begin{lemma}\label{lem:base-estimates}
For all sufficiently large $N$,
\begin{equation}\label{eq:kappa-derivative}
 0\le\kappa'\le C N^{-6}\kappa^2.
\end{equation}
Let \(I\) be an interval that intersects at most one transition
interval, possibly extending into the adjacent intervals where
the frequency is constant.
If it meets the $j$-th
transition, set $k=k_j$; if it contains only one constant frequency,
let $k$ be that frequency.
Then
\begin{equation}\label{eq:base-estimates}
 \norm{B_0}_{\A(I)}\le Ck^2,\qquad
 \norm{\tr B_0}_{\A(I)}\le Ck^2N^{-2}.
\end{equation}
Moreover,
\begin{equation}\label{eq:flatness}
 \kappa(t)\to\infty,\qquad a(t)\kappa(t)^p\to0\quad(p>0)
 \quad\text{as }t\uparrow0,
\end{equation}
and
\begin{equation}\label{eq:a-small}
 \sup_{-1\le t<0}a(t)\kappa(t)^2\le N^{-2}.
\end{equation}
The function $g$ is a smooth increasing bijection from $(-\infty,0)$
onto $\R$.
\end{lemma}

\begin{proof}
      Differentiating $\kappa$ gives
\[
 \kappa'=(k_{j+1}-k_j)\frac{k_j}{N^2}
 S'\!\left(\frac{k_j(t-t_j)}{N^2}\right),
 \qquad 0\le k_{j+1}-k_j\le3N^{-4}k_j.
\]
This proves \eqref{eq:kappa-derivative}, since on the constant intervals
the derivative is zero and \eqref{eq:ratio} makes $k_j$ and $\kappa$ comparable.
On the $j$-th transition,
\[
 |\theta_j'|\le Ck_jN^{-2},\qquad
 |\theta_j''|\le Ck_j^2N^{-4}.
\]
For convenience, denote $c=\cos \theta_j$ or $\sin \theta_j$, and $d=k_j$ or $k_{j+1}$.
After differentiation, we have 
\begin{equation}\label{eq:c-derivatives}
 |c|\le1,\qquad |c'|\le CkN^{-2},\qquad
 |c''|\le Ck^2N^{-4}.
\end{equation}
for either amplitude and on the adjacent constant intervals.
Since $d$ equals $k_j$ or $k_{j+1}$,
\begin{equation}\label{eq:frequency-nearby}
 k\le \kappa,d\le(1+3N^{-4})k,
 \qquad |\kappa-d|\le CkN^{-4}.
\end{equation}

For the term $c(t)\cos(dx)$, the nonzero entries of
$a^{-1}D^2(ac\cos(dx))$ are
\begin{equation}\label{eq:one-wave-matrix}
 \begin{split}
 tt:&\quad [c''-2\kappa c'+(\kappa^2-\kappa')c]\cos(dx),\\
 tx=xt:&\quad d(\kappa c-c')\sin(dx),\\
 xx:&\quad -d^2c\cos(dx).
 \end{split}
\end{equation}
The nonzero entries of
$a^{-1}D^2(ac\cos(dy))$ are similar to \eqref{eq:one-wave-matrix} with $x$ replaced by $y$.
% The $\A(I)$ norms of $\cos(dx)$ and $\sin(dx)$ are both one.
Each displayed coefficient is bounded by $Ck^2$ using
\eqref{eq:c-derivatives}--\eqref{eq:frequency-nearby}. Summing all
entries proves the first bound in
\eqref{eq:base-estimates}.
The trace coefficient for either $a^{-1}D^2(ac\cos(dx))$ or $a^{-1}D^2(ac\cos(dy))$ is
\[
 c''-2\kappa c'+(\kappa^2-d^2)c-\kappa'c.
\]
Its four terms have absolute values at most
\[
 Ck^2N^{-4},\quad Ck^2N^{-2},\quad
 |\kappa-d|(\kappa+d)\le Ck^2N^{-4},\quad Ck^2N^{-6},
\]
respectively. Their sum is at most $Ck^2N^{-2}$, which proves the trace bound.

On $[t_j,t_{j+1}]$, $k_j\le\kappa\le k_{j+1}$, so
$\kappa(t)\to\infty$ as $t\uparrow0$. Further,
\begin{equation}\label{eq:integal}
    g(t_{j+1})-g(t_j)=\int_{t_j}^{t_{j+1}}\kappa(s)\,ds
 \ge k_j\frac{2N^2}{k_j}=2N^2.
\end{equation}
For $t\in[t_j,t_{j+1}]$, this gives
\[
 a(t)\kappa(t)^p
 \le N^{-2}k_1^{-2}\e^{-g(t_1)}
        \e^{-2N^2(j-1)}(N^4+j+1)^{2p}.
\]
For fixed $N,p$, the right side therefore tends to zero, proving
\eqref{eq:flatness}. \eqref{eq:integal} gives $g(t)\to
\infty$ as $t\uparrow0$. For $t<t_1$, $g'=k_1$, so
$g(t)\to-\infty$ as $t\to-\infty$. Finally, $g'=\kappa>0$.
These facts prove the asserted bijection and its smooth inverse.
By \eqref{eq:kappa-derivative},
\[
 \frac{d}{dt}\log(a\kappa^2)
 =-\kappa+\frac{2\kappa'}{\kappa}
 \le-\kappa(1-2CN^{-6})\le-\frac\kappa2
\]
once $2CN^{-6}\le1/2$. Thus $a\kappa^2$ decreases. Since
$-1<t_1$, $g(-1)=0$ and $\kappa(-1)=k_1$, its value at $-1$
is $N^{-2}$. This proves
\eqref{eq:a-small}.
\end{proof}

\section{A right inverse for the frozen operator}
% \label{sec:constant}\label{sec:frozen}

The following lemma is a key ingredient in the proof. 
In this section, we construct a right inverse for a linear equation with parameters frozen in $t$. 
In the next section, we apply this construction on short time intervals and control the resulting errors 
to solve the linear equation with variable coefficients.

For a positive integer $d$, define
\begin{equation*}\label{eq:T-matrices}
 T_x(d)=d^2\begin{pmatrix}
 \cos(dx)&\sin(dx)&0\\
 \sin(dx)&-\cos(dx)&0\\0&0&0
 \end{pmatrix},\quad
 T_y(d)=d^2\begin{pmatrix}
 \cos(dy)&0&\sin(dy)\\
 0&0&0\\\sin(dy)&0&-\cos(dy)
 \end{pmatrix}.
\end{equation*}
If the amplitude of $\cos(dx)$ is constant and
$\kappa=d$, its contribution to $B_0$ in
\eqref{eq:one-wave-matrix} is precisely that amplitude times
$T_x(d)$. Both matrices have trace zero.

\begin{lemma}\label{lem:constant-inverse}
Let $d_1,d_2$ be positive integers and $k=\min(d_1,d_2)$. Suppose
\[
 \max(d_1,d_2)\le11k/10,\qquad k\le\kappa_0\le11k/10.
\]
For $c_1,c_2\ge0$ with $c_1^2+c_2^2=1$, set
\begin{equation*}\label{eq:constant-operator}
 \begin{split}
 \overline B&=c_1 T_x(d_1)+c_2 T_y(d_2),\\
 \overline P(\psi,z)&=
 ((\partial_t-\kappa_0)^2+\Delta_{x,y})z-\overline B:D^2\psi.
 \end{split}
\end{equation*}
There is a real linear right inverse $R_0$, such that, for every real $r \in \A(\R)$,
$\overline P(R_0r)=r$ with
\begin{equation}\label{eq:constant-bound}
 \norm\psi_{2,k;\R}+k^{-2}\norm z_{2,k;\R}
 \le Ck^{-2}\norm r_{\A(\R)},\qquad (\psi,z)=R_0r.
\end{equation}
The constant is independent of $k,d_1,d_2,\kappa_0,c_1,c_2$.
\end{lemma}

\begin{proof}
\textbf{Step 1: Choose frequencies and derive the coefficient equations.}
For now, assume $c_1,c_2\ge\eta>0$. The number $\eta$ will
be fixed at the end of the proof. Set
\[
 E_k=\{m\in\Z^2:|m|\le5k/4\}.
\]
In this proof, $e_1=(1,0),e_2=(0,1)$ denote the fixed coordinate
vectors. Given $m\in E_k\setminus\{0\}$, choose $j\in\{1,2\}$
so that $d_j|m_j|$ is largest, taking $j=1$ in case of a tie. Define
\begin{equation*}\label{eq:frequency-choice}
 \sigma=\operatorname{sgn}(m_j),\qquad q=m+\sigma d_je_j.
\end{equation*}
We will choose $\psi_q$ to produce the required coefficient $r_m$.

The sine and cosine in the $j$-th matrix shift the frequency
$q$ to $q\pm d_je_j$. One is $m$; the other is $m+2\sigma d_je_j$.
The latter lies outside $E_k$, since $|m_j|+2d_j\ge2k>5k/4$.
The other matrix, with index $h\ne j$, produces $q\pm d_he_h$.
Since $d_j|m_j|\ge d_h|m_h|$, we have
\begin{align*} \label{eq:other-frequencies}
 |m+\sigma d_je_j\pm d_he_h|^2
 &=|m|^2+d_j^2+d_h^2+2d_j|m_j|\pm2d_hm_h\notag\\
 &\ge |m|^2+d_j^2+d_h^2
 \ge2k^2>25k^2/16.
\end{align*}
Thus, among the frequencies in $E_k$, $\psi_q$ contributes only at $m$.

We show that the map $m\mapsto q(m)$ is injective. By the choice of $j$,
\[
\begin{aligned}
 j=1:\quad d_1|q_1|-d_2|q_2|
 &=d_1|m_1|-d_2|m_2|+d_1^2\ge d_1^2,\\
 j=2:\quad d_1|q_1|-d_2|q_2|
 &=d_1|m_1|-d_2|m_2|-d_2^2\le-d_2^2.
\end{aligned}
\]
These two cases are disjoint, so $q$ determines $j$.
Since $\sigma=\operatorname{sgn}(q_j)$, we can then recover
$m=q-\sigma d_je_j$. This proves injectivity.
Replacing $m$ by $-m$ leads to
\begin{equation}\label{eq:conjugate-choice}
 q(-m)=-q(m).
\end{equation}

To compute the coefficient equation, it suffices to write the $x$ direction explicitly. Since
\[
 T_x(d_1):D^2\psi
 =d_1^2\{\cos(d_1x)(\psi_{tt}-\psi_{xx})
                           +2\sin(d_1x)\psi_{tx}\},
\]
$\psi_q\e^{\ii(q_1x+q_2y)}$ contributes, at the
frequency $q-\sigma d_1 e_1$,
\[
 -c d_1^2\left\{\frac12(\psi_q''+q_1^2\psi_q)
                +2\frac{\ii\sigma}{2}\,\ii q_1\psi_q'\right\}
=-\frac{cd_1^2}{2}
       (\psi_q''-2\sigma q_1\psi_q'+q_1^2\psi_q).
\]
For $y$, replace $q_1$ by $q_2$ and
$d_1$ by $d_2$. So the coefficient at the desired frequency
is exactly
\begin{equation}\label{eq:coefficient-equation}
 -\frac{c_jd_j^2}{2}(\partial_t-|q_j|)^2\psi_q.
\end{equation}
\textbf{Step 2: Construct and estimate \(\psi\) in low frequencies.}
For $d>0$ and bounded measurable $f$, define
\begin{equation}\label{eq:G-definition}
 (G_df)(t)=\int_0^\infty s e^{-ds}f(t+s)\,ds.
\end{equation}
Thus $G_df=L_d*f$, where
$L_d(s)=-s e^{ds}\mathbf{1}_{\{s<0\}}$.
Since $L_d$ is continuous and $L_d'$ has a unit jump at zero,
\( D^2L_d=(-2d-d^2s)e^{ds}\mathbf{1}_{\{s<0\}}+\delta_0.\)
A direct calculation gives $(D-d)^2L_d=\delta_0$, hence
$(\partial_t-d)^2G_df=f$ in distributions. Moreover,
\[
 \|L_d\|_{L^1}=d^{-2},\qquad
 \|L_d'\|_{L^1}\le 2d^{-1},\qquad
 |D^2L_d|(\mathbb R)\le 4.
\]
Using $(G_df)'=L_d'*f$ and $(G_df)''=(D^2L_d)*f$,
we obtain
\begin{equation*}\label{eq:G-bounds}
 d^2\|G_df\|_\infty
 +d\|(G_df)'\|_\infty
 +\|(G_df)''\|_\infty
 \le C\|f\|_\infty,
\end{equation*}
with an absolute constant $C\ge 7$. The derivatives of \(G_df\) are understood in the weak sense;
no derivative of $f$ is required.

For $m\ne0$ in $E_k$, use its chosen $j,q$ and set
\begin{equation}\label{eq:psi-nonzero}
 \psi_q=-\frac{2}{c_jd_j^2}G_{|q_j|}r_m.
\end{equation}
\eqref{eq:coefficient-equation} and
\eqref{eq:G-definition} show that its contribution at $m$ is
$r_m$. Step~1 proves that none of the other chosen Fourier coefficients
contributes at this target.

For $m=0$, choose $j_*$ with $c_{j_*}\ge1/\sqrt2$. Set
\begin{equation}\label{eq:psi-zero}
 \psi_{d_{j_*}e_{j_*}}=\psi_{-d_{j_*}e_{j_*}}
       =-\frac1{c_{j_*}d_{j_*}^2}G_{d_{j_*}}r_0.
\end{equation}
Each of these two Fourier coefficients produces $r_0/2$ at zero by
\eqref{eq:coefficient-equation}, so their sum gives $r_0$.
Their other frequencies are $\pm2d_{j_*}e_{j_*}$ and
$\pm d_{j_*}e_{j_*}\pm d_he_h$, where $h\ne j_*$. The norms are at least $\sqrt{2}k>\frac{5}{4}k$.
All unspecified Fourier coefficients of $\psi$ are set to zero.

For each chosen coefficient,
\begin{equation*}\label{eq:q-sizes}
 k\le|q_j|\le|q|\le|m|+d_j\le\frac{47}{20}k.
\end{equation*}
By \eqref{eq:psi-nonzero} and \eqref{eq:psi-zero}, the bounds
for $\psi_q,\psi_q',\psi_q''$ are, respectively,
\[
 2\eta^{-1}k^{-4}\norm{r_m}_\infty,\qquad
 4\eta^{-1}k^{-3}\norm{r_m}_\infty,\qquad
 8\eta^{-1}k^{-2}\norm{r_m}_\infty.
\]
Multiplication by $k^2+|q|^2$, $k+|q|$, and $1$ respectively, shows that its
contribution to $\norm\psi_{2,k;\R}$ is at most
$C\eta^{-1}k^{-2}\norm{r_m}_\infty$. Injectivity allows summation, so
\begin{equation}\label{eq:psi-estimate}
 \norm\psi_{2,k;\R}
 \le C\eta^{-1}k^{-2}\sum_{m\in E_k}\norm{r_m}_\infty.
\end{equation}
\textbf{Step 3: Construct and estimate $z$ in high frequencies.}
Set
\begin{equation}\label{eq:remaining-data}
 r^*=r+\overline B:D^2\psi.
\end{equation}
Steps~1--2 imply $r_m^*=0$ on $E_k$. The $\A(\mathbb{R})$-norm sum of entries in $T_x(d)$
is $4d^2$, and the same holds for $T_y(d)$.
Thus $\norm{\overline B}_{\A(\mathbb{R})}\le Ck^2$. Using
\eqref{eq:derivative-control} and \eqref{eq:psi-estimate} gives
\begin{equation}\label{eq:rstar-bound}
 \norm{r^*}_{\A(\R)}
 \le\norm r_{\A(\R)}+Ck^2\norm\psi_{2,k;\R}
 \le C(1+\eta^{-1})\norm r_{\A(\R)}.
\end{equation}

For $|m|>5k/4$, define
\begin{equation*}\label{eq:K-definition}
 K_{|m|,\kappa_0}(s)=-\frac{1}{2|m|}
 \begin{cases}
 e^{(\kappa_0-|m|)s},&s\ge0,\\
 e^{(\kappa_0+|m|)s},&s<0.
 \end{cases}
\end{equation*}
As in Step 2, differentiation in the sense of distributions gives
\begin{equation}\label{eq:K-equation}
 ((D-\kappa_0)^2-|m|^2)K_{|m|,\kappa_0}=\delta_0.
\end{equation}
Direct integration gives
\begin{equation}\label{eq:K-integrals}
 \|K\|_{L^1}=\frac{1}{|m|^2-\kappa_0^2},
 \qquad
 \|K'\|_{L^1}=\frac{1}{|m|},
 \qquad
 |D^2K|(\mathbb R)=2.
\end{equation}
Set
\[
 z_m=K_{|m|,\kappa_0}*r_m^*\quad(m\notin E_k),
 \qquad z_m=0\quad(m\in E_k).
\]
For $m\notin E_k$, we have $\kappa_0/|m|\le22/25$ and
$k\le4|m|/5$. Together with \eqref{eq:K-integrals},
we obtain
\[
 (k^2+|m|^2)\|z_m\|_\infty
 +(k+|m|)\|z_m'\|_\infty+\|z_m''\|_\infty
 \le C\|r_m^*\|_\infty.
\]
Summing over $m$ and using \eqref{eq:rstar-bound} gives
\[
 \|z\|_{2,k;\mathbb R}
 \le C\|r^*\|_{\mathcal A(\mathbb R)}
 \le C(1+\eta^{-1})\|r\|_{\mathcal A(\mathbb R)}.
\]
By \eqref{eq:K-equation} and the summable derivative bounds,
$((\partial_t-\kappa_0)^2+\Delta_{x,y})z=r^*$ in distributions.
Together with \eqref{eq:remaining-data}, this proves
$\overline P(\psi,z)=r$ and \eqref{eq:constant-bound},
with a constant depending on $\eta$ currently.
Equations \eqref{eq:conjugate-choice} and \eqref{eq:psi-zero},
together with the real kernels and coefficients, preserve conjugate
symmetry. Hence $\psi,r^*,z$ are real.

\noindent \textbf{Step 4: Construct an inverse when one amplitude is zero.}
Suppose only j-th matrix is present, with amplitude $c_j$ satisfying $1/\sqrt2 \le c_j \le 1$.
The construction below gives a right inverse with an absolute bound, uniform throughout this interval.
Keep the same $k,\kappa_0,E_k$ as above.
For every nonzero $m\in E_k$, set
\[
 \sigma=\begin{cases}
 \operatorname{sgn}(m_j),&m_j\ne0,\\
 \operatorname{sgn}(m_h),&m_j=0,
 \end{cases}\qquad q=m+\sigma d_je_j,
 \quad h\ne j.
\]
In both cases,
$|q_j|=|m_j|+d_j$ and $\operatorname{sgn}(q_j)=\sigma$.
Thus $q$ determines $m=q-\operatorname{sgn}(q_j)d_je_j$;
the map $m \to q(m)$ is injective and $q(m)=-q(-m)$.
The other frequency is $m+2\sigma d_je_j$,
with $|m_j|+2d_j\ge2k$ outside $E_k$.
For $m=0$, define $\psi_{d_{j}e_{j}}=\psi_{-d_{j}e_{j}}$ as in \eqref{eq:psi-zero}, 
which gives the same contribution after frequency shifts. 
Other frequencies $\pm 2d_je_j$ are outside $E_k$.

As in Step 2, using \eqref{eq:psi-nonzero}--\eqref{eq:psi-zero}, for $m \in E_k$, we have
\[
 (k^2+|q|^2)\norm{\psi_q}_\infty
       +(k+|q|)\norm{\psi_q'}_\infty+\norm{\psi_q''}_\infty
\le \frac{C}{c_jk^2}\norm{r_m}_\infty.
\]
Injectivity and $c_j^{-1}\le\sqrt2$ therefore give
\[
 \norm\psi_{2,k;\R}
 \le \frac{C}{c_jk^2}\sum_{m\in E_k}\norm{r_m}_\infty
 \le Ck^{-2}\norm r_{\A(\R)}.
\]
For high frequencies, take $j=1$; the case $j=2$ follows
by interchanging $x,y$. Define
\[
 r^*=r+c_1T_x(d_1):D^2\psi.
\]
As in Step 3, we know that $r_m^*=0$ on $E_k$, while
\[
 \norm{r^*}_{\A(\R)}
 \le\norm r_{\A(\R)}+4c_1d_1^2\norm{D^2\psi}_{\A(\R)}
 \le C\norm r_{\A(\R)}.
\]
Set $z_m=K_{|m|,\kappa_0}*r_m^*$ off $E_k$ and $z_m=0$ on
$E_k$. Equations \eqref{eq:K-equation}--\eqref{eq:K-integrals}
then give
\[
 ((\partial_t-\kappa_0)^2+\Delta_{x,y})z
                  -c_1T_x(d_1):D^2\psi=r,\qquad
 \norm\psi_{2,k;\R}+k^{-2}\norm z_{2,k;\R}
                  \le Ck^{-2}\norm r_{\A(\R)}.
\]
Thus $R_*r=(\psi,z)$ is a right inverse with an absolute bound.

\noindent\textbf{Step 5: Treat small nonzero amplitude by convergent series.}
If $c_2<\eta<1/\sqrt2$, then $c_1>1/\sqrt2$. Define $P_*$ by
\[
 \overline P=P_*+D,\qquad D(\psi,z)=-c_2 T_y(d_2):D^2\psi.
\]
So $P_*$ has one amplitude zero. Step~4 then gives
\[
 k^{-2}\norm{D R_*r}_{\A}
 \le k^{-2}c_2\,4d_2^2\norm{D^2(R_*r)_\psi}_{\A}
 \le Cc_2\,k^{-2}\norm r_{\A}.
\]
Here $(R_*r)_\psi$ denotes the first component. Fix $\eta$ so that $\eta<1/\sqrt2$ and
$C\eta\le1/2$. The operator series
\[
 (I+D R_*)^{-1}=\sum_{\nu=0}^\infty(-D R_*)^\nu
\]
then converges and has norm at most two. It follows that
\begin{equation*}\label{eq:small-amplitude-inverse}
 R_0=R_*\sum_{\nu=0}^\infty(-D R_*)^\nu
\end{equation*}
is a right inverse, since
\[
 \overline P R_0
 =(P_*R_*+DR_*)(I+DR_*)^{-1}
 =(I+DR_*)(I+DR_*)^{-1}=I.
\]
If $c_1<\eta$, interchange $x$ and $y$. Since
$\eta$ now is fixed, a single constant works in
\eqref{eq:constant-bound}.
\end{proof}

\section{Patching local inverses for the variable-coefficient operator}
\label{sec:global}
After setting $w=a(h_0+z)$ and dividing \eqref{eq:exact} by $a$, the cofactor term has linear part
\begin{equation*}\label{eq:P}
P(\psi,z)=\tr B_z+(\tr D^2\psi)(\tr B_0)-D^2\psi:B_0.
\end{equation*}
Its coefficients vary with $t$. We freeze amplitudes and $\kappa$ on short intervals, patch local solutions together, and estimate coefficient differences and cutoff terms.
A Neumann series then gives a uniform right inverse for the contraction argument in Section~\ref{sec:nonlinear}.

Choose smooth cutoffs $\chi,\zeta:\R\to[0,1]$ such that
\[
 \chi=1\text{ on }[-1,1],\quad\supp\chi\subset(-2,2),\qquad
 \zeta=1\text{ on }[-2,2],\quad\supp\zeta\subset(-3,3).
\]
For $i\in\Z$, set
\begin{equation*}\label{eq:cover}
 \begin{split}
 s_i&=g^{-1}(Ni),\qquad V_i=\{t<0:|g(t)/N-i|<3\},\\
 \chi_i(t)&=\frac{\chi(g(t)/N-i)}{\sum_{j\in\Z}\chi(g(t)/N-j)},
 \qquad \zeta_i(t)=\zeta(g(t)/N-i).
 \end{split}
\end{equation*}
Thus,
\begin{equation}\label{eq:partition-properties}
 0\le\chi_i,\zeta_i\le1,\qquad
 \sum_i\chi_i=1,\qquad \zeta_i\chi_i=\chi_i.
\end{equation}
Their supports are compact subsets of $V_i$. Since $g$ is an
increasing bijection, each $V_i$ has compact closure in $(-\infty,0)$.
If $V_i\cap V_j\ne\varnothing$, then $|i-j|<6$; hence each
$V_i$ intersects at most ten other intervals $V_j$.
Between two consecutive transition intervals, $\kappa$ is constant. The $g$-length
of this constant-frequency interval is at least $N^2$, since
\[
\int_{t_j+N^2/k_j}^{t_{j+1}}\kappa(t)\,dt
=N^2\frac{k_{j+1}}{k_j}\ge N^2.
\]
For $N>6$, we have $|g(V_i)|=6N<N^2$, so $V_i$ intersects
at most one transition interval. Thus, on $V_i$, $h_0$ has
the following form with fixed positive integers $d_1,d_2$,
\begin{equation}\label{eq:local-h0}
 h_0=c_1(t)\cos(d_1x)+c_2(t)\cos(d_2y),\qquad
 c_1,c_2\ge0,\,c_1^2+c_2^2=1.
\end{equation}
If there is a transition, $d_1,d_2$ are its two frequencies,
assigned to their actual coordinate directions. If there is only one
frequency, assign the same frequency to the other direction with zero
amplitude. 

Set $K_i:=\min(d_1,d_2)$ and
$\kappa_i=\kappa(s_i)$.
We use $K_i$ to distinguish the local scale on $V_i$ from the
frequency sequence $k_j$ in Section~\ref{sec:algebra}.
If $V_i$ meets the $j$-th transition, then
\( \{d_1,d_2\}=\{k_j,k_{j+1}\}\), and \(K_i=k_j.\)
If $V_i$ lies in a constant part with frequency $k_j$, we take \(d_1=d_2=k_j\), and \(K_i=k_j.\)
Thus each local scale $K_i$ is a term of the frequency sequence
from Section~\ref{sec:algebra}, but its index in that sequence need not be $i$.

By \eqref{eq:ratio}, for large $N$,
\begin{equation}\label{eq:local-comparable}
 K_i\le\kappa(t),d_1,d_2\le11K_i/10\quad(t\in V_i).
\end{equation}
In particular, if $V_i\cap V_j\ne\varnothing$, choosing a point of their
intersection in \eqref{eq:local-comparable} gives
\begin{equation}\label{eq:scale-overlap}
 \frac{10}{11}\le\frac{K_i}{K_j}\le\frac{11}{10}.
\end{equation}
Also $6N=\int_{V_i}\kappa\ge K_i|V_i|$, so
\begin{equation}\label{eq:interval-length}
 |V_i|\le6N/K_i.
\end{equation}
The cutoffs also satisfy
\begin{equation}\label{eq:cutoff-derivatives}
 \norm{\chi_i'}_\infty+\norm{\zeta_i'}_\infty\le CK_i/N,
 \quad
 \norm{\chi_i''}_\infty+\norm{\zeta_i''}_\infty\le CK_i^2/N^2.
\end{equation}

For real Fourier series on $(-\infty,0)\times\T^2$, define
\begin{equation}\label{eq:global-norms1}
 \norm{(\psi,z)}_\XX=\sup_i\left(\norm\psi_{2,K_i;V_i}
                         +K_i^{-2}\norm z_{2,K_i;V_i}\right),
\end{equation}
\begin{equation}\label{eq:global-norms2}
     \norm r_\YY=\sup_i K_i^{-2}\norm r_{\A(V_i)}.
\end{equation}
Let $\XX$ be the space of pairs $(\psi,z)$ with finite \eqref{eq:global-norms1} norm,
and let $\YY$ be the space of function $r$ with finite \eqref{eq:global-norms2} norm.
Both spaces are complete, whose proof is similar to that of the spaces with norm \eqref{eq:A-norm} and \eqref{eq:second-norm}.

Write $p=(\psi,z)$. By \eqref{eq:A-norm}, \eqref{eq:second-norm}, \eqref{eq:derivative-control}, and \eqref{eq:global-norms1}, we have
\begin{equation}\label{eq:local-controls}
 \begin{array}{lll}
 \norm\psi_\A\le K_i^{-2}\norm p_\XX,&
 \norm{D\psi}_\A\le CK_i^{-1}\norm p_\XX,&
 \norm{D^2\psi}_\A\le C\norm p_\XX,\\
 \norm z_\A\le\norm p_\XX,&
 \norm{Dz}_\A\le CK_i\norm p_\XX,&
 \norm{D^2z}_\A\le CK_i^2\norm p_\XX,
 \end{array}
\end{equation}
where $\norm{\cdot}_\A$ means $\norm{\cdot}_{\A(V_i)}$. From
\eqref{eq:B-entries}, $\kappa\le CK_i$ and
$|\kappa'|\le CK_i^2$, it follows that
\begin{equation}\label{eq:Bz-control}
 \norm{B_z}_{\A(V_i)}\le CK_i^2\norm p_\XX.
\end{equation}
\eqref{eq:algebra} and \eqref{eq:base-estimates} now yield
\[
 K_i^{-2}\norm{Pp}_{\A(V_i)}
 \le K_i^{-2}\norm{B_z}_\A
       +CK_i^{-2}\norm{D^2\psi}_\A\norm{B_0}_\A
 \le C\norm p_\XX.
\]
Thus $P:\XX\to\YY$ is bounded.

We first compare $P$ with $\overline P_i$ defined below and estimate the terms
produced by cutoffs. On each $V_i$, set
\[
 \begin{split}
 \overline B_i&=c_1(s_i)T_x(d_1)+c_2(s_i)T_y(d_2),\\
 \overline P_i(\psi,z)
 &=((\partial_t-\kappa_i)^2+\Delta_{x,y})z
                         -\overline B_i:D^2\psi.
 \end{split}
\]

\begin{lemma}\label{lem:local-errors}
The following estimates hold for $N$ as above.
\begin{equation}\label{eq:B-error}
 \norm{B_0-\overline B_i}_{\A(V_i)}\le CK_i^2/N.
\end{equation}
For every real pair $p=(\psi,z)$ on $\R\times\T^2$ with
$\norm\psi_{2,K_i;\R}+K_i^{-2}\norm z_{2,K_i;\R}<\infty$,
\begin{equation}\label{eq:P-error-estimate}
 K_i^{-2}\norm{(P-\overline P_i)p}_{\A(V_i)}
 \le\frac CN\left(\norm\psi_{2,K_i;\R}
                         +K_i^{-2}\norm z_{2,K_i;\R}\right),
\end{equation}
and
\begin{equation}\label{eq:commutator-estimate}
 K_i^{-2}\norm{[P,\zeta_i]p}_{\A(V_i)}
 \le\frac CN\left(\norm\psi_{2,K_i;\R}
                         +K_i^{-2}\norm z_{2,K_i;\R}\right),
\end{equation}
where
\[
 [P,\zeta_i]p=P(\zeta_i\psi,\zeta_i z)-\zeta_iP(\psi,z).
\]
The constant $C$ is independent of $i$ and $N$.
\end{lemma}

\begin{proof}
For either amplitude $c=c_1$ or $c_2$, 
\eqref{eq:c-derivatives} and \eqref{eq:interval-length} give
\[
 |c(t)-c(s_i)|\le |V_i|\sup_{V_i}|c'|
 \le(6N/K_i)(CK_i/N^2)\le C/N.
\]
For the $x$ term, its contribution to
$B_0-\overline B_i$ has coefficients
\begin{equation}\label{eq:B-error-entries}
 \begin{split}
 tt:&\ c_1''-2\kappa c_1'-\kappa'c_1
                +(\kappa^2-d_1^2)c_1+d_1^2(c_1-c_1(s_i)),\\
 tx=xt:&\ d_1\{(\kappa-d_1)c_1-c_1'+d_1(c_1-c_1(s_i))\},\\
 xx:&\ -d_1^2(c_1-c_1(s_i)).
 \end{split}
\end{equation}
They multiply $\cos(d_1x),\sin(d_1x),\cos(d_1x)$, respectively.
The $y$ term has exactly the displayed coefficients with
$d_1,c_1,x$ replaced by $d_2,c_2,y$ and entries $tx,xx$ replaced
by $ty,yy$. In the first line of \eqref{eq:B-error-entries}, the
bounds for five terms are
\[
 CK_i^2N^{-4},\quad CK_i^2N^{-2},\quad CK_i^2N^{-6},
 \quad CK_i^2N^{-4},\quad CK_i^2N^{-1},
\]
respectively. In the second line, they are bounded by $CK_i^2N^{-4},CK_i^2N^{-2}$ and
$CK_i^2N^{-1}$, respectively; the third line is bounded by $CK_i^2N^{-1}$.
Summing proves \eqref{eq:B-error}.

Since $\kappa$ and $\kappa_i$ lie between the same adjacent
frequencies, \eqref{eq:ratio}, \eqref{eq:kappa-derivative}, and
\eqref{eq:local-comparable} give
\[
 |\kappa-\kappa_i|\le CK_iN^{-4},\qquad
 |\kappa^2-\kappa_i^2-\kappa'|\le CK_i^2N^{-4}.
\]
Subtracting the two operators gives
\begin{equation*}\label{eq:P-error-formula}
 (P-\overline P_i)(\psi,z)=-2(\kappa-\kappa_i)z_t
       +(\kappa^2-\kappa_i^2-\kappa')z+(\tr B_0)(\tr D^2\psi)
               -(B_0-\overline B_i):D^2\psi.
\end{equation*}
Use \eqref{eq:derivative-control} on a pair $(\psi, z)$,
after division by $K_i^2$, the first two terms have norm at most
\[
 CK_i^{-2}\left(
 K_iN^{-4}\frac{\norm z_{2,K_i;\R}}{K_i}
 +K_i^2N^{-4}\frac{\norm z_{2,K_i;\R}}{K_i^2}
 \right)
 \le CN^{-4}K_i^{-2}\norm z_{2,K_i;\R}.
\]
The remaining terms, by \eqref{eq:base-estimates} and
\eqref{eq:B-error}, have divided norm at most
\[
 CK_i^{-2}(K_i^2N^{-2}+K_i^2N^{-1})
                    \norm\psi_{2,K_i;\R}
 \le CN^{-1}\norm\psi_{2,K_i;\R}.
\]
Together these prove \eqref{eq:P-error-estimate}.

To estimate the commutator, let $e_t=(1,0,0)$. By the product rule,
\[
 D^2(\zeta_i\psi)-\zeta_iD^2\psi
 =\zeta_i'(e_t\otimes D\psi+D\psi\otimes e_t)
       +\zeta_i''\psi\,e_t\otimes e_t.
\]
Substitution into the definition of $P$ gives
\begin{equation*}\label{eq:commutator-formula}
 [P,\zeta_i]p
 =2\zeta_i'z_t+(\zeta_i''-2\kappa\zeta_i')z
 +((\tr B_0)I_3-B_0):
 \bigl[\zeta_i'(e_t\otimes D\psi+D\psi\otimes e_t)
                 +\zeta_i''\psi\,e_t\otimes e_t\bigr].
\end{equation*}
The first two terms, after division by $K_i^2$, are bounded by
\[
 CK_i^{-2}\left[
 \frac{K_i}{N}\frac{\norm z_{2,K_i;\R}}{K_i}
 +\left(\frac{K_i^2}{N^2}+\frac{K_i^2}{N}\right)
             \frac{\norm z_{2,K_i;\R}}{K_i^2}
 \right]
 \le\frac CN K_i^{-2}\norm z_{2,K_i;\R}.
\]
By \eqref{eq:base-estimates},
\[
 \norm{(\tr B_0)I_3-B_0}_{\A(V_i)}
 \le 3\norm{\tr B_0}_{\A(V_i)}
                  +\norm{B_0}_{\A(V_i)}
 \le CK_i^2.
\]
The remaining terms, after division by $K_i^2$, therefore have norm at most
\[
 C\left[
 \frac{K_i}{N}\frac{\norm\psi_{2,K_i;\R}}{K_i}
 +\frac{K_i^2}{N^2}\frac{\norm\psi_{2,K_i;\R}}{K_i^2}
 \right]
 \le\frac CN\norm\psi_{2,K_i;\R}.
\]
This proves \eqref{eq:commutator-estimate}.
\end{proof}

\begin{proposition}\label{prop:global-inverse}
For all sufficiently large $N$, there is a real linear operator
$\mathcal R:\YY\to\XX$ satisfying
\begin{equation*}\label{eq:global-inverse}
 P\mathcal R=I_\YY,\qquad
 \norm{\mathcal R}_{\YY\to\XX}\le C,
\end{equation*}
with $C$ independent of $N$.
\end{proposition}

\begin{proof}
Let $R_i$ be the right inverse of $\overline P_i$ from Lemma~\ref{lem:constant-inverse}. Extend $\chi_i r$ by zero
outside $V_i$ on $\R$. Since $0\le\chi_i\le1$,
\[
 K_i^{-2}\norm{\chi_i r}_{\A(\R)}
 \le K_i^{-2}\norm r_{\A(V_i)}
 \le\norm r_\YY.
\]
Set $p_i=R_i(\chi_i r)$ and define
\begin{equation*}\label{eq:approximate-inverse}
 \mathcal S r=\sum_{i\in\Z}\zeta_i p_i.
\end{equation*}
The sum is locally finite. \eqref{eq:constant-bound} gives
\[
 \norm{(p_i)_\psi}_{2,K_i;\R}
 +K_i^{-2}\norm{(p_i)_z}_{2,K_i;\R}
 \le C\norm r_\YY.
\]
By \eqref{eq:cutoff-norm} and \eqref{eq:cutoff-derivatives},
the same bound holds with $p_i$ replaced by $\zeta_i p_i$.
For $V_i\cap V_j\ne\varnothing$, \eqref{eq:scale-overlap}
and the norm definition \eqref{eq:second-norm} then give
\[
 \norm{\zeta_i(p_i)_\psi}_{2,K_j;V_j}
 +K_j^{-2}\norm{\zeta_i(p_i)_z}_{2,K_j;V_j}
 \le C\norm r_\YY.
\]
Since $\zeta_i p_i$ is supported in $V_i$, at most eleven terms
contribute on each $V_j$. The triangle inequality and
\eqref{eq:global-norms1} therefore yield
\begin{equation}\label{eq:S-bound}
 \norm{\mathcal S r}_\XX\le C\norm r_\YY.
\end{equation}

Using $\overline P_i p_i=\chi_i r$ and
$\zeta_i\chi_i=\chi_i$, we obtain
\begin{equation}\label{eq:error-identity}
 P\mathcal S r=r+\mathcal E r,\qquad
 \mathcal E r=\sum_i\left[
       \zeta_i(P-\overline P_i)p_i+[P,\zeta_i]p_i
       \right].
\end{equation}
Denote
$e_i=\zeta_i(P-\overline P_i)p_i+[P,\zeta_i]p_i$.
Lemma~\ref{lem:local-errors}, together with
\eqref{eq:partition-properties} and \eqref{eq:constant-bound} gives
\[
 K_i^{-2}\norm{e_i}_{\A(V_i)}
 \le CN^{-1}\norm r_\YY.
\]
On a fixed $V_j$, the triangle inequality gives
\[
 \begin{split}
 K_j^{-2}\sum_m\norm{(\mathcal E r)_m}_{\infty,V_j}
 &\le\sum_{i:V_i\cap V_j\ne\varnothing}
     \left(\frac{K_i}{K_j}\right)^2
           K_i^{-2}\sum_m\norm{(e_i)_m}_{\infty,V_i}
 \le\frac CN\norm r_\YY.
 \end{split}
\]
Taking the supremum over $j$ proves
\begin{equation}\label{eq:E-bound}
 \norm{\mathcal E}_{\YY\to\YY}\le C/N.
\end{equation}

Choose $N$ large enough that the right side of
\eqref{eq:E-bound} is at most $1/2$. The Neumann series for
$(I+\mathcal E)^{-1}$ therefore converges in operator norm on
the Banach space $\YY$ and has norm at most two. Define
\begin{equation*}\label{eq:R-definition}
 \mathcal R=\mathcal S\sum_{\nu=0}^\infty(-\mathcal E)^\nu.
\end{equation*}
Due to the boundedness of $P$, \eqref{eq:error-identity} and
\eqref{eq:S-bound} then give
\[
 \begin{aligned}
 P\mathcal R
 &=(I+\mathcal E)\sum_{\nu=0}^\infty(-\mathcal E)^\nu
 =I_\YY,\\
 \norm{\mathcal R}_{\YY\to\XX}
 &\le\norm{\mathcal S}_{\YY\to\XX}
           \sum_{\nu=0}^\infty
                    \norm{\mathcal E}_{\YY\to\YY}^\nu
 \le2\norm{\mathcal S}_{\YY\to\XX}\le C.
 \end{aligned}
\]
All cutoffs, coefficient matrices, and local inverses preserve real
functions, so does $\mathcal R$.
\end{proof}

\begin{remark}\label{rem}
For fixed \(N\), we suppress the dependence of the spaces \(\mathcal X,\mathcal Y\) and the operators \(P,\mathcal S,\mathcal E,\mathcal R\) on \(N\). 
More precisely, there exist
constants $N_0,C_S,C_E$, depending only on the fixed functions $S,\chi,\zeta$,
such that for every integer $N\ge N_0$ the operators $\mathcal S_N$ and
$\mathcal E_N$ satisfy
\[
P_N\mathcal S_N=I_{\mathcal Y_N}+\mathcal E_N,
\qquad
\|\mathcal S_N\|\le C_S,
\qquad
\|\mathcal E_N\|\le C_E/N.
\]
Thus, for every $N\ge\max\{N_0,2C_E\}$, the operator
\[
\mathcal R_N
=
\mathcal S_N(I_{\mathcal Y_N}+\mathcal E_N)^{-1}
\]
is a right inverse of $P_N$ and satisfies $\|\mathcal R_N\|\le 2C_S$.
Whenever $N$ is increased below, all $N$-dependent objects are understood
to be constructed with the new value of $N$; all constants in the
estimates remain unchanged.
\end{remark}

% The formulas allow data to affect other time intervals. The decay
% needed at $t=0$ follows from the weights in $\XX$, not from a
% support claim for these integral operators.

\section{Solving the exact nonlinear equation}\label{sec:nonlinear}

Section~\ref{sec:global} provides a uniformly bounded right inverse for \(P\), 
the linear part of the normalized cofactor expression in \eqref{eq:exact}. 
\eqref{eq:exact} also contains the trace, nonlinear cofactor, and the determinant. 
We estimate these remaining terms and apply contraction map to construct a real pair \((\psi,z)\in\mathcal X\) with \(\|(\psi,z)\|_{\mathcal X}\) small,
such that \(\psi\) and \(w=a(h_0+z)\) satisfy \eqref{eq:exact} almost everywhere on \(-1<t<0\).

The construction of the inverse in Section~\ref{sec:global} uses all $t<0$, but the functions
in \Cref{thm:main} are only needed for $-1<t<0$. Since $a$ grows as
$t\to-\infty$, set
\begin{equation*}\label{eq:tau}
 \tau(t)=\mathbf1_{(-1,0)}(t)a(t)^2.
\end{equation*}
It will only multiply a determinant and will never be differentiated.
Therefore a bounded measurable cutoff suffices. By
\eqref{eq:a-small} and \eqref{eq:local-comparable},
\begin{equation}\label{eq:tau-bound}
 \sup_{t<0}\tau(t)\kappa(t)^4\le N^{-4},\qquad
 \norm\tau_{L^\infty(V_i)}\le N^{-4}K_i^{-4}.
\end{equation}

Define
\begin{equation}\label{eq:F}
 F(\psi,z)=\cof(I_3+D^2\psi):(B_0+B_z)
                          +\tau\det(B_0+B_z).
\end{equation}
We aim to solve $F(\psi,z)=0$ on $t<0$. On $-1<t<0$, since
$D^2w=a(B_0+B_z)$, homogeneity of the three-dimensional determinant
gives
\begin{equation*}\label{eq:F-exact}
 \cof(I_3+D^2\psi):D^2w+\det D^2w
 =a\cof(I_3+D^2\psi):(B_0+B_z)
 +a^3\det(B_0+B_z)
 =a F(\psi,z).
\end{equation*}
Thus a solution of \eqref{eq:F} gives the required equation \eqref{eq:exact}.
Using \eqref{eq:cofactor}, the cofactor part expands as
\[
 \begin{split}
 \cof(I_3+D^2\psi):(B_0+B_z)
 &=\tr B_0+\tr B_z+((\tr D^2\psi)I_3-D^2\psi):B_0\\
 &\quad+((\tr D^2\psi)I_3-D^2\psi):B_z +\cof D^2\psi:(B_0+B_z).
 \end{split}
\]
Separating the constant term, we have
\begin{equation}\label{eq:F-expansion}
 F(p)=F(0)+Pp+\mathcal N(p),\qquad
 F(0)=\tr B_0+\tau\det B_0,
\end{equation}
where
\begin{equation}\label{eq:N}
 \mathcal N(p)=((\tr D^2\psi)I_3-D^2\psi):B_z+\cof D^2\psi:(B_0+B_z)+\tau\{\det(B_0+B_z)-\det B_0\}.
\end{equation}
In particular, $\mathcal N(0)=0$.
\begin{lemma}\label{lem:nonlinear-estimates}
There are constants $C_0,C_1$ independent of large $N$ such that
\begin{equation}\label{eq:F0-bound}
 \norm{F(0)}_\YY\le C_0N^{-2}.
\end{equation}
For $0<r\le1$ and $\norm p_\XX,\norm q_\XX\le r$,
\begin{equation}\label{eq:N-lipschitz}
 \norm{\mathcal N(p)-\mathcal N(q)}_\YY
 \le C_1(r+N^{-4})\norm{p-q}_\XX.
\end{equation}
\end{lemma}

\begin{proof}
On a fixed $V_i$, abbreviate the norm $\A(V_i)$ by $\A$.
\eqref{eq:base-estimates}, \eqref{eq:polynomial-size},
and \eqref{eq:tau-bound} give
\[
 \begin{split}
 \norm{F(0)}_\A
 &\le\norm{\tr B_0}_\A
            +\norm\tau_{\infty,V_i}\norm{\det B_0}_\A \le CK_i^2N^{-2}+N^{-4}K_i^{-4}\,C(K_i^2)^3\\
 &=CK_i^2(N^{-2}+N^{-4}).
 \end{split}
\]
Divide by $K_i^2$, use $N^{-4}\le N^{-2}$, and take the
supremum over $i$. This proves \eqref{eq:F0-bound}.

Take two pairs $p=(\psi,z)$ and $q=(\widetilde\psi,\widetilde z)$.
\eqref{eq:local-controls} and \eqref{eq:Bz-control}
give
\begin{equation}\label{eq:nonlinear-sizes}
\begin{gathered}
 \norm{D^2\psi}_\A+\norm{D^2\widetilde\psi}_\A
 \le Cr,
 \qquad
 \norm{D^2\psi-D^2\widetilde\psi}_\A
 \le C\norm{p-q}_\XX,
 \\[4pt]
 \norm{B_z}_\A+\norm{B_{\widetilde z}}_\A
 \le CK_i^2r,
 \qquad
 \norm{B_z-B_{\widetilde z}}_\A
 \le CK_i^2\norm{p-q}_\XX,
 \\[4pt]
 \norm{B_0+B_z}_\A+\norm{B_0+B_{\widetilde z}}_\A
 \le CK_i^2(1+r).
\end{gathered}
\end{equation}
We now apply \eqref{eq:nonlinear-sizes} to each term of \eqref{eq:N}.
$\mathcal N(p)-\mathcal N(q)$ consists of three parts. First,
\[
 \begin{split}
 \mathcal{E}_1:= &((\tr D^2\psi)I_3-D^2\psi):B_z
       -((\tr D^2\widetilde\psi)I_3-D^2\widetilde\psi):B_{\widetilde z}\\
 =&((\tr(D^2\psi-D^2\widetilde\psi))I_3-(D^2\psi-D^2\widetilde\psi)):B_z
 +((\tr D^2\widetilde\psi)I_3-D^2\widetilde\psi):
                                  (B_z-B_{\widetilde z}).
 \end{split}
\]
Since
$\norm{(\tr D^2\psi)I_3-D^2\psi}_\A\le4\norm{D^2\psi}_\A$,
$\mathcal{E}_1$ is therefore bounded by
\[
 C\norm{p-q}_\XX\,K_i^2r
       +Cr\,K_i^2\norm{p-q}_\XX
 \le CK_i^2r\norm{p-q}_\XX.
\]
Second,
\[
 \begin{split}
 \mathcal{E}_2:=&\cof D^2\psi:(B_0+B_z)-\cof D^2\widetilde\psi:(B_0+B_{\widetilde z})\\
 =&(\cof D^2\psi-\cof D^2\widetilde\psi):(B_0+B_z)
                  +\cof D^2\widetilde\psi:(B_z-B_{\widetilde z}).
 \end{split}
\]
From \eqref{eq:cofactor-difference} and
\eqref{eq:nonlinear-sizes},
\[
 \norm{\cof D^2\psi-\cof D^2\widetilde\psi}_\A
      \le Cr\norm{p-q}_\XX,
 \qquad \norm{\cof D^2\widetilde\psi}_\A\le Cr^2.
\]
Consequently, $\mathcal{E}_2$ is bounded by
\[
 Cr\norm{p-q}_\XX\,K_i^2(1+r)
       +Cr^2K_i^2\norm{p-q}_\XX
 \le CK_i^2(r+r^2)\norm{p-q}_\XX.
\]
By \eqref{eq:det-difference} and \eqref{eq:nonlinear-sizes}, the last term $\mathcal{E}_3:=\tau\{\det(B_0+B_z)-\det(B_0+B_{\widetilde z})\}$ is  also bounded,
\[
 \begin{split}
 &\norm{\tau\{\det(B_0+B_z)-\det(B_0+B_{\widetilde z})\}}_\A\\
 \le &C\norm\tau_{\infty,V_i}
       \bigl(\norm{B_0+B_z}_\A+\norm{B_0+B_{\widetilde z}}_\A\bigr)^2
       \norm{B_z-B_{\widetilde z}}_\A\\
 \le &C(N^{-4}K_i^{-4})\,K_i^4(1+r)^2
                                  K_i^2\norm{p-q}_\XX
 = CN^{-4}K_i^2(1+r)^2\norm{p-q}_\XX.
 \end{split}
\]
Adding $\mathcal{E}_1$, $\mathcal{E}_2$, $\mathcal{E}_3$ together,
dividing by $K_i^2$, and taking the supremum over $i$ gives
\[
 \norm{\mathcal N(p)-\mathcal N(q)}_\YY
 \le C\{r+r^2+N^{-4}(1+r)^2\}\norm{p-q}_\XX.
\]
Since $r\le1$, $r^2\le r$ and $(1+r)^2\le4$, we prove \eqref{eq:N-lipschitz}.
\end{proof}

\begin{proposition}\label{prop:nonlinear-solution}
For all sufficiently large $N$, there is a real pair
$p=(\psi,z)\in\XX$ with
\begin{equation}\label{eq:solution-size}
 F(\psi,z)=0\quad\text{a.e. on }t<0,
 \qquad\norm p_\XX\le CN^{-2}.
\end{equation}
\end{proposition}

\begin{proof}
Let $C_R$ be a bound for $\norm{\mathcal R}_{\YY\to\XX}$ (we can take $C_R=2C_S$, the uniform bound from \Cref{rem}) and set
\begin{equation*}\label{eq:radius}
 r_N=2C_RC_0N^{-2}.
\end{equation*}
The constants are independent of $N$. Increase $N$ until
\begin{equation*}\label{eq:contraction-conditions}
 r_N\le1,\qquad C_RC_1(r_N+N^{-4})\le\frac14.
\end{equation*}
For $\norm p_\XX\le r_N$ define
\begin{equation*}\label{eq:T-map}
 T(p)=-\mathcal R\bigl(F(0)+\mathcal N(p)\bigr).
\end{equation*}
Using $\mathcal N(0)=0$, \Cref{lem:nonlinear-estimates} gives
\[
 \norm{T(p)}_\XX \le C_RC_0N^{-2}+C_RC_1(r_N+N^{-4})\norm p_\XX \le 3r_N/4,\]
\[\norm{T(p)-T(q)}_\XX\le\tfrac14\norm{p-q}_\XX.\] 
Thus $T$ is a contraction mapping.
Define a closed ball $B^*=\{p\in\XX:\norm p_\XX\le r_N\}$, which is complete.
By Banach contraction theorem,
$T$ therefore has a unique fixed point in $B^*$. In particular,
starting with $p^{(0)}=0$, the iteration
\begin{equation}\label{eq:iteration}
 p^{(j+1)}=-\mathcal R\bigl(F(0)+\mathcal N(p^{(j)})\bigr)
\end{equation}
converges in $\XX$ to $p$ with $\norm p_\XX\le r_N$.
The contraction estimate also gives
\[
 \norm{p^{(j+1)}-p^{(j)}}_\XX
 \le4^{-j}\norm{p^{(1)}}_\XX
 \le\tfrac12r_N4^{-j}.
\]
Applying $P\mathcal R=I$ to \eqref{eq:iteration} gives, for $j\ge1$,
\[
 Pp^{(j)}=-F(0)-\mathcal N(p^{(j-1)}).
\]
Thus, \eqref{eq:F-expansion} and \eqref{eq:N-lipschitz} imply
\[
 F(p^{(j)})=\mathcal N(p^{(j)})-\mathcal N(p^{(j-1)}),\]
\[ \norm{F(p^{(j)})}_\YY
 \le C_1(r_N+N^{-4})\,\tfrac12r_N4^{1-j}\longrightarrow0.
\]
On $B^*$, \eqref{eq:F-expansion} also yields
\[
 \norm{F(p^{(j)})-F(p)}_\YY
 \le\bigl(\norm P_{\XX\to\YY}+C_1(r_N+N^{-4})\bigr)
                       \norm{p^{(j)}-p}_\XX\longrightarrow0.
\]
Consequently, $F(p)=0$ in $\YY$. For each $i$, its Fourier
coefficients vanish almost everywhere on $V_i$, so
$F(p)=0$ almost everywhere on $V_i\times\T^2$.
The countable cover $\{V_i\}$ gives $F(p)=0$ almost everywhere on $(-\infty, 0)\times\T^2$.
All iterations and their limit are real. This proves
\eqref{eq:solution-size}.
\end{proof}

\section{Proof of \Cref{thm:main}}\label{sec:completion}

We first extend the functions across $t=0$ and prove
Proposition~\ref{prop:three}. We then prove Theorem~\ref{thm:main}
in dimension three and extend to higher dimensions.
Take the real pair $(\psi,z)$ supplied by
Proposition~\ref{prop:nonlinear-solution}. Extend $\psi$ by zero for
$t\ge0$, and define on $B_1^3$
\begin{equation}\label{eq:extension}
 \phi=Q_3+\psi,\qquad
 w(t,x,y)=
 \begin{cases}
  a(t)\bigl(h_0(t,x,y)+z(t,x,y)\bigr),&t<0,\\
  0,&t\ge0.
 \end{cases}
\end{equation}

\begin{lemma}\label{lem:extension}
For all sufficiently large $N$, the functions in
\eqref{eq:extension} belong to $C^{1,1}(B_1^3)$ and satisfy
\begin{equation*}\label{eq:final-small}
 \norm{D^2\phi-I_3}_{L^\infty(B_1^3)}
 +\norm{D^2w}_{L^\infty(B_1^3)}\le CN^{-2},
\end{equation*}
where $C$ is independent of $N$. Moreover, $\phi=Q_3$ and $w=0$
on $B_{1,+}^3$, while $w(t,0,0)>0$ for $-1<t<0$.
% The function $w$ also takes negative values in $B_1^3$.
\end{lemma}

\begin{proof}
By \eqref{eq:local-comparable}, \eqref{eq:local-controls}, and
\eqref{eq:solution-size}, the following estimates hold for almost every \((t,x,y)\in (V_i\cap(-1,0))\times\mathbb T^2\).
\begin{equation}\label{eq:psi-pointwise}
\begin{gathered}
|\psi(t,x,y)|\le CN^{-2}\kappa(t)^{-2},\;
|D\psi(t,x,y)|\le CN^{-2}\kappa(t)^{-1},\;
|D^2\psi(t,x,y)|\le CN^{-2},\\[3pt]
|z(t,x,y)|\le CN^{-2},\qquad
|Dz(t,x,y)|\le CN^{-2}\kappa(t).
\end{gathered}
\end{equation}

For the two nonnegative amplitudes in $h_0$,
$c_1^2+c_2^2=1$, and hence
$c_1+c_2\le\sqrt2$. Thus, by \eqref{eq:h0}, \eqref{eq:local-h0} and \eqref{eq:c-derivatives},
\[
 |h_0|\le\sqrt2,
 \qquad |\partial_xh_0|+|\partial_yh_0|\le C\kappa,
 \qquad |\partial_th_0|\le C\kappa N^{-2}.
\]
The product rule, with $a'=-\kappa a$, now gives
\begin{gather*}
 |w|=a|h_0+z|\le Ca,\qquad |w_t|\le a\{\kappa|h_0+z|+|(h_0)_t|+|z_t|\}
                                                \le Ca\kappa,\\
 |w_x|+|w_y|\le a\{|(h_0)_x|+|(h_0)_y|+|z_x|+|z_y|\}
                                                \le Ca\kappa.
\end{gather*}
Finally, on $-1<t<0$, we have $D^2w=a(B_0+B_z)$, so almost everywhere
\begin{equation}\label{eq:w-pointwise}
 |w|\le Ca,\qquad |Dw|\le Ca\kappa,\qquad
 |D^2w|\le Ca\kappa^2\le CN^{-2}.
\end{equation}
For the Hessian estimate, we use
$\norm{B_0}_{\A(V_i)}\le CK_i^2$,
$\norm{B_z}_{\A(V_i)}\le CK_i^2N^{-2}$, and
\eqref{eq:a-small}.

\eqref{eq:flatness}, \eqref{eq:psi-pointwise}, and
\eqref{eq:w-pointwise} show that $\psi, w$ and $D\psi, Dw$ tend uniformly to zero as $t\uparrow0$.
For $F=\psi$ or $w$, \eqref{eq:psi-pointwise} and
\eqref{eq:w-pointwise} give $|D^2F|\le L:=CN^{-2}$ almost everywhere
in the negative half-ball. The $W^{1,\infty}$ estimate in \cite[Lemma~7.8(iii)]{Tartar2007},
applied to $e\cdot DF$ for each unit vector $e$, shows that $DF$
is $L$-Lipschitz, thus
$|DF(t,x,y)-DF(s,x,y)|\le L(s-t)$ for $t<s<0$.
Taking $s\uparrow0$ and using the gradient limit above gives
$|DF(t,x,y)|\le L|t|$. Since $F(0,x,y)=0$, integration in $t$ yields
\begin{equation*}\label{eq:quadratic-interface}
 |F(t,x,y)|\le \int_t^0 L|s|\,ds=\frac{L}{2}t^2.
\end{equation*}
For any $X_0\in B_1^3\cap\{t=0\}$, this gives
$|F(X)-F(X_0)|\le \frac{L}{2}|X-X_0|^2$ on the negative side.
The same bound holds on the positive side since $F=0$ there.
Thus $DF(X_0)=0$.
If $X=(t,x,y)$ lies in the negative half-ball and
$Y\in B_1^3\cap\{t\ge0\}$, then
\[
 |DF(X)-DF(Y)|=|DF(X)|\le L|t|\le L|X-Y|.
\]
Hence $DF$ is $L$-Lipschitz on $B_1^3$. This proves
$\psi,w\in C^{1,1}(B_1^3)$ and
\begin{equation*}
 \norm{D^2\psi}_{L^\infty(B_1^3)}
 +\norm{D^2w}_{L^\infty(B_1^3)}\le CN^{-2}
\end{equation*}
for a fixed constant $C$ independent of $N$.

At $(x,y)=(0,0)$, $h_0(t)=\cos\theta_j(t)+\sin\theta_j(t)$. In a
transition, its value is $\cos\theta_j(t)+\sin\theta_j(t)\ge1$. 
On a constant interval, the value is 1.
Choose $N$ large enough that $|z|\le1/2$ by
\eqref{eq:psi-pointwise}. Since $a(t)>0$,
\(w(t,0,0)=a(t)(h_0(t,0,0)+z(t,0,0))\ge a(t)/2>0\),
when \(-1<t<0\).
% The function $w$ also takes a negative value. Because $-1/2<t_1$, at
% $X=(-1/2,\pi/k_1,0)$ we have $h_0(X)=\cos\pi=-1$ and
% $h_0(X)+z(X)\le-1/2$. This point lies in the ball for large $N$,
% since $|X|^2=1/4+\pi^2/k_1^2<1$. Hence $w(X)<0$.
\end{proof}

\begin{proof}[Proof of Proposition~\ref{prop:three}]
Given $\delta>0$, choose $N$ sufficiently large that
$CN^{-2}\le\delta$ in Lemma~\ref{lem:extension}.
On $B_1^3\cap\{t<0\}$ we have
\[
 D^2\phi=I_3+D^2\psi,\qquad
 D^2w=a(B_0+B_z),\qquad \tau=a^2.
\]
Consequently, by \eqref{eq:F} and Proposition~\ref{prop:nonlinear-solution},
\[
 \begin{aligned}
 \cof(D^2\phi):D^2w+\det D^2w&=a\bigl\{\cof(I_3+D^2\psi):(B_0+B_z)
                  +a^2\det(B_0+B_z)\bigr\}\\
 &=aF(\psi,z)=0
 \qquad\text{almost everywhere}.
 \end{aligned}
\]
On $B_{1,+}^3$, the same identity follows from $w=0$.
Since $B_1^3\cap\{t=0\}$ has Lebesgue measure zero,
\eqref{eq:boxed} holds almost everywhere on $B_1^3$.
The regularity, the bound \eqref{eq:small-hessians}, and the assertions
about the values of $\phi,w$ follow from Lemma~\ref{lem:extension}.
\end{proof}

\begin{proof}[Proof of Theorem~\ref{thm:main} in dimension three]
Fix $0<\varepsilon<1$. Apply Proposition~\ref{prop:three} with
$\delta=\varepsilon/4$ and set
\[
 u_3=\phi+w,\qquad v_3=\phi-w.
\]
Both functions belong to $C^{1,1}(B_1^3)$.
For almost every $X\in B_1^3$ and every $\xi\in\R^3$,
\[
 |\xi^T(D^2\phi-I_3\pm D^2w)\xi|
 \le(|D^2\phi-I_3|+|D^2w|)|\xi|^2
 \le\frac\varepsilon4|\xi|^2.
\]
Consequently,
\begin{equation}\label{eq:three-spectral}
 (1-\varepsilon/4)I_3\le D^2u_3,D^2v_3
                      \le(1+\varepsilon/4)I_3
 \quad\text{almost everywhere},
\end{equation}
implying both $u$ and $v$ are uniformly convex on $B_1^3$.

By \eqref{eq:odd-det} and
\eqref{eq:boxed}, almost everywhere on $B_1^3$, we have
\begin{equation}\label{eq:equal-determinants}
 \begin{aligned}
 \det D^2u_3-\det D^2v_3
 &=\det(D^2\phi+D^2w)-\det(D^2\phi-D^2w)\\
 &=2\bigl\{\cof(D^2\phi):D^2w+\det D^2w\bigr\}=0.
 \end{aligned}
\end{equation}

Choose a null set $E_3\subset B_1^3$ such that both Hessians
exist and \eqref{eq:three-spectral}, \eqref{eq:equal-determinants}
hold at every point of $B_1^3\setminus E_3$. Define
\begin{equation*}\label{eq:density-three}
 f_3(X)=
 \begin{cases}
  \det D^2u_3(X), & X\in B_1^3\setminus E_3,\\
  1,              & X\in E_3.
 \end{cases}
\end{equation*}
Then $f_3$ is measurable and agrees almost everywhere with both
$\det D^2u_3$ and $\det D^2v_3$. \eqref{eq:three-spectral} gives
\[
 1-\varepsilon<(1-\varepsilon/4)^3\le f_3\le(1+\varepsilon/4)^3<1+\varepsilon.
\]

Applying \eqref{eq:area} to $u_3$ and $v_3$,
\begin{equation*}\label{eq:measure-equality-three}
 \begin{aligned}
 M_{u_3}(E)&=|Du_3(E)|=\int_E\det D^2u_3(X)\,dX=\int_E f_3(X)\,dX\\
 &=\int_E\det D^2v_3(X)\,dX
  =|Dv_3(E)|=M_{v_3}(E)
 \end{aligned}
\end{equation*}
for every Borel set $E\subset B_1^3$. In particular, both measures
assign zero mass to $B_1^3\cap\{t=0\}$. Finally,
Proposition~\ref{prop:three} gives $u_3=v_3=Q_3$ on $B_{1,+}^3$, and
$u_3(t,0,0)-v_3(t,0,0)=2w(t,0,0)>0$ for $-1<t<0$.
This proves all four assertions in dimension three.
\end{proof}

\begin{proof}[Proof of Theorem~\ref{thm:main} in dimensions $n>3$]
Let $u_3,v_3,f_3$ be the functions just constructed for the same
$\varepsilon$.
Write $X=(X',X'')\in\R^3\times\R^{n-3}$ and define
\begin{equation*}\label{eq:lift}
 u_n(X)=u_3(X')+\tfrac12|X''|^2,\qquad
 v_n(X)=v_3(X')+\tfrac12|X''|^2.
\end{equation*}
If $X\in B_1^n$, then $|X'|\le|X|<1$. Their gradients are
$(Du_3(X'),X'')$ and $(Dv_3(X'),X'')$, which are Lipschitz and injective.
Their Hessians are
\[
 D^2u_n=\begin{pmatrix}D^2u_3&0\\0&I_{n-3}\end{pmatrix},
 \qquad
 D^2v_n=\begin{pmatrix}D^2v_3&0\\0&I_{n-3}\end{pmatrix}, \quad \text{almost everywhere}.
\]
Set $E_n= B_1^n \cap E_3\times [-1,1]^{n-3}$. Fubini theorem gives
$|E_n| \le 2^{n-3}|E_3|=0$. For every \(X=(X',X'')\in B_1^n\setminus E_n\), \eqref{eq:equal-determinants}
holds at \(X'\).
 Thus, $\det(D^2u_n)=\det(D^2v_n)$ holds almost everywhere in $B_1^n$.
\eqref{eq:three-spectral} implies
\[
 (1-\varepsilon/4)I_n\le D^2u_n,D^2v_n\le(1+\varepsilon/4)I_n
 \quad\text{almost everywhere}.
\]
Set $f_n(X)=f_3(X')$. For almost every $X\in B_1^n$, we have
\[
 \det D^2u_n(X)=\det D^2u_3(X')=f_n(X)=\det D^2v_n(X).
\]
Its bounds are still $1-\varepsilon$ and $1+\varepsilon$.
Applying \eqref{eq:area} for every Borel $E\subset B_1^n$,
\[
 M_{u_n}(E)=\int_E\det D^2u_n(X)\,dX
   =\int_E f_3(X')\,dX =\int_E f_n(X)\,dX=M_{v_n}(E).
\]
This proves \eqref{eq:main-measure}, including for sets intersecting $\{X_1=0\}$.

If $X_1>0$, 
\(u_n(X)=v_n(X)=Q_n(X),\) and $f_n(X)=1$.
For $-1<t<0$, \(u_n(t,0,\ldots,0)-v_n(t,0,\ldots,0)=u_3(t,0,0)-v_3(t,0,0)=2w(t,0,0)>0.\)
All assertions of Theorem~\ref{thm:main} now hold in every dimension $n>3$.
\end{proof}

\begin{proof}[Proof of \Cref{cor:interval}]
Set
\[
    \varepsilon=\frac{M-m}{M+m},
    \qquad
    s=\left(\frac{m+M}{2}\right)^{1/n},
\]
and let $u_n,v_n,f_n$ be given by Theorem~\ref{thm:main}
for this choice of $\varepsilon$.

Fix $X_0\in\Omega$ and choose $0<r<R$ such that
\(B_r(X_0)\subset\Omega\subset B_R(X_0).\)
Define $T(X)=(X-X_0)/R$ on $B_R(X_0)$ and
\[
    u(X)=sR^2 u_n(T(X)),\qquad
    v(X)=sR^2 v_n(T(X)),\qquad
    f(X)=s^n f_n(T(X)).
\]
The restrictions of $u$ and $v$ to $\Omega$ are uniformly convex
and have Lipschitz gradients. Since $Du(X)=sR\,Du_n(T(X))$ and $Dv(X)=sR\,Dv_n(T(X))$, every Borel set $E\subset\Omega$ satisfies
\[
    M_u(E)
    =(sR)^n M_{u_n}(T(E))
    =\int_E f(X)\,dX
    =(sR)^n M_{v_n}(T(E))
    =M_v(E).
\]
The density bounds follow from
$s^n(1-\varepsilon)=m$ and $s^n(1+\varepsilon)=M$.

Let $e_1=(1,0,\ldots,0)$. By part~\textup{(2)} of
Theorem~\ref{thm:main}, the two solutions agree on the nonempty
open set \(\{X\in\Omega:(X-X_0)\cdot e_1>0\}\),
while $u(X_0-he_1)>v(X_0-he_1)$ for every $0<h<r$.
Hence $u$ and $v$ are distinct on $\Omega$.
\end{proof}

\section*{Acknowledgement}
The AI tool ChatGPT (OpenAI) was used in the preparation of this manuscript. The author
has reviewed all AI-assisted content and takes full responsibility for the final manuscript.

\bibliographystyle{amsplain}
\bibliography{references}

\end{document}